\documentclass[10pt]{amsart}

\usepackage{amssymb,amsmath,amsthm,amsfonts,longtable}
\usepackage[all]{xy}

\theoremstyle{plain}
\newtheorem{Theorem}{Theorem}
\newtheorem{Proposition}[Theorem]{Proposition}
\newtheorem{Lemma}[Theorem]{Lemma}
\newtheorem{Question}[Theorem]{Question}
\newtheorem{Corollary}[Theorem]{Corollary}
\newtheorem{Conjecture}[Theorem]{Conjecture}
\newtheorem{Problem}[Theorem]{Problem}

\theoremstyle{definition}
\newtheorem{Remark}[Theorem]{Remark}
\newtheorem{Observation}[Theorem]{Observation}

\begin{document}
\title{On infinite groups generated by two quaternions}
\author{Diego Rattaggi}
\thanks{Supported by the Swiss National Science Foundation, No.\ PP002--68627}
%\address{Dorfstrasse 38, 
%CH--6005 Luzern, Switzerland}
\email{rattaggi@hotmail.com}
\date{\today}
\begin{abstract}
Let $x$, $y$ be two integral quaternions of norm $p$ and $l$,
respectively, where $p$, $l$ are distinct odd prime numbers.
We investigate the structure of $\langle x,y \rangle$,
the multiplicative group generated by $x$ and $y$.
Under a certain condition which excludes $\langle x,y \rangle$
from being free or abelian, we show for example that $\langle x,y \rangle$,
its center, commutator subgroup and abelianization are finitely
presented infinite groups.
We give many examples where our condition is satisfied
and compute as an illustration a finite presentation of the group
$\langle 1+j+k, 1+2j \rangle$
having these two generators and seven relations.
In a second part, we study the basic question whether there exist 
commuting quaternions $x$ and $y$ for fixed $p$, $l$, using results
on prime numbers of the form $r^2 + m s^2$ and a simple invariant for commutativity.
\end{abstract}
\maketitle
\setcounter{section}{-1}
\section{Introduction} \label{intro}
Let $p$, $l$ be two distinct odd prime numbers and
$x$, $y$ two integral Hamilton quaternions whose norms are in the set
$\{ p^r l^s : r, s \in \mathbb{N}_0 \} \setminus \{ 1 \}$.
We are interested in the structure of the multiplicative group generated
by $x$ and $y$. These groups $\langle x,y \rangle$ are always infinite since
$x$ and $y$ have infinite order by the assumption made on the norms.

We first consider the case where 
$x$ and $y$ do not commute. It is well-known that certain pairs $x$, $y$
generate a free group of rank two, e.g.\ $\langle 1+2i, 1+2k \rangle \cong F_2$.
However, if the norms of $x$ and $y$ are not powers of the same prime
number $p$, then the structure of $\langle x,y \rangle$ is unknown in general.
Nevertheless, we have shown in \cite[Proposition~27]{Rattaggi2} that 
the group $\langle 1+2i, 1+4k \rangle$ is \emph{not} free
by establishing an explicit relation of length $106$.
This holds in a more general situation.
We will describe in Section~\ref{SectionPrel} a homomorphism $\psi_{p,l}$
defined on the group of invertible rational quaternions,
and give the definition of a group $\Gamma_{p,l}$ 
such that $\langle \psi_{p,l}(x), \psi_{p,l}(y) \rangle$
is a subgroup of $\Gamma_{p,l}$.
It follows from results in \cite{Rattaggi2} that if $x$, $y$
satisfy some technical condition on the parity of their coefficients
and if the group $\langle \psi_{p,l}(x), \psi_{p,l}(y) \rangle$ has finite index
in $\Gamma_{p,l}$, then $\langle x,y \rangle$ is not free.
Here we will use similar techniques to get more precise information on some groups $\langle x,y \rangle$
and naturally related subgroups or quotients.
For instance, under the assumptions mentioned above which imply that $\langle x,y \rangle$ is not free,
we will prove in Theorem~\ref{Thm11} that $\langle x,y \rangle$ as well as
its center, its commutator subgroup and its abelianization are
finitely presented infinite groups. This contrasts to the case where $\langle x,y \rangle \cong F_2$,
since then the center is finite (in fact trivial) and the commutator subgroup
is not finitely presented (in fact not finitely generated). 
It also contrast to the case where $x$ and $y$ commute, since then
the commutator subgroup of $\langle x,y \rangle$ is trivial.
However, we will show in Theorem~\ref{Thm11} that 
our groups $\langle x,y \rangle$
contain non-abelian free subgroups $F_2$ of infinite index and 
free abelian subgroups $\mathbb{Z} \times \mathbb{Z}$ of infinite index.

Our constructions will be illustrated in Section~\ref{SectionEx}
for the concrete example $x = 1+j+k$ of norm $p=3$ and $y = 1+2j$ of
norm $l=5$. 
In particular, we will compute a finite presentation of 
$\langle 1+j+k,1+2j  \rangle$ having seven relations,
and determine its center, which turns out to be 
$\langle 3^4, 5^4 \rangle = \langle 81, 625 \rangle < \mathbb{Q}^{\ast}$, a group isomorphic to
$\mathbb{Z} \times \mathbb{Z}$.
We guess that the finiteness assumption made on the index of $\langle \psi_{p,l}(x), \psi_{p,l}(y) \rangle$
in $\Gamma_{p,l}$ is a rather restrictive condition in general. Nevertheless
we are able to give many explicit examples where it holds
(at least for small $p$ and $l$).
Table~\ref{Table1} in Section~\ref{lists} describes 191 selected examples of such pairs
$x$, $y$ of norm $p$ and $l$, respectively, for 56 distinct pairs $(p,l)$
satisfying $3 \leq p < l < 100$.

If $x$ and $y$ commute, then the structure of $\langle x,y \rangle$
is less challenging. We get an abelian group like
\[
\langle 1+2i, 1+4i \rangle \cong \mathbb{Z} \times \mathbb{Z}
\]
or
\[
\langle 1+2i, -1-2i \rangle = \langle 1+2i, -1 \rangle \cong \mathbb{Z} \times \mathbb{Z}_2,
\]
where we always use the notation $\mathbb{Z}_n := \mathbb{Z}/n\mathbb{Z}$
for the cyclic group of order $n$.
However, we are interested to know for which pairs $p$, $l$ there
exist commuting integral quaternions $x$, $y$ of norm $p$ and $l$ at all.
We will study this problem in Section~\ref{SectionComm}
and give some partial general answers using congruence conditions for $p$ and $l$.

\section{Preliminaries} \label{SectionPrel}
Throughout the whole article, let $p$, $l$ be any pair of distinct odd prime numbers.
The goal of this section is to define and describe the family of groups
$\Gamma_{p,l}$ mentioned in the introduction.
These groups are closely related to certain finitely generated multiplicative subgroups
of invertible rational Hamilton quaternions.

For a commutative ring $R$ with unit, let
\[
\mathbb{H}(R) = 
\{ x_0 + x_1 i + x_2 j + x_3 k : x_0, x_1, x_2, x_3 \in R \}
\]
be the ring of Hamilton quaternions over $R$, i.e.\ $1,i,j,k$ is a free basis,
and the multiplication is determined by the identities $i^2 = j^2 = k^2 = -1$ and
$ij = -ji =k$.
Let $\overline{x} := x_0 - x_1 i - x_2 j - x_3 k$ be the \emph{conjugate}
of $x =  x_0 + x_1 i + x_2 j + x_3 k \in \mathbb{H}(R)$, and
\[
|x|^2 := x \overline{x} = \overline{x} x = x_0^2 + x_1^2 + x_2^2 + x_3^2 \in R 
\]
its \emph{norm}. 
Note that $|xy|^2 = |x|^2 |y|^2$ for all $x, y \in \mathbb{H}(R)$.
We denote by $\Re(x) := x_0$ the ``real part'' of $x$.
Let $R^{\ast}$ be the multiplicative group of invertible elements in the ring $R$.
We will mainly use the two groups $\mathbb{H}(\mathbb{Q})^{\ast} = \mathbb{H}(\mathbb{Q}) \setminus \{ 0 \}$
and $\mathbb{Q}^{\ast} = \mathbb{Q} \setminus \{ 0 \}$.
Let
\[
\mathbb{H}(R)_1 := \{ x \in \mathbb{H}(R)^{\ast} : |x|^2 = 1  \}
\]
be the subgroup of quaternions of norm $1$. 

If $K$ is a field, let as usual 
$\mathrm{PGL}_2(K) = \mathrm{GL}_2(K)/Z\mathrm{GL}_2(K)$ be
the quotient of the group of invertible $(2 \times 2)$-matrices over $K$
by its center.
We write brackets $[A] \in \mathrm{PGL}_2(K)$ to denote the image of the matrix 
$A \in \mathrm{GL}_2(K)$
under the quotient homomorphism $\mathrm{GL}_2(K) \to \mathrm{PGL}_2(K)$.
Let $\mathbb{Q}_p$, $\mathbb{Q}_l$ be the field of $p$-adic and $l$-adic
numbers, respectively, and fix elements 
$c_p, d_p \in \mathbb{Q}_p$ and $c_l, d_l \in \mathbb{Q}_l$
such that 
\[
c_p^2 + d_p^2 + 1 = 0 \in \mathbb{Q}_p
\text{ \, and \, } 
c_l^2 + d_l^2 + 1 = 0 \in \mathbb{Q}_l.
\]
For $q \in \{p,l\}$ let $\psi_q$ be the homomorphism of groups
$\mathbb{H}(\mathbb{Q})^{\ast} \to \mathrm{PGL}_2(\mathbb{Q}_q)$
defined by 
\[
\psi_q(x) := \left[
\begin{pmatrix}
x_0 + x_1 c_q + x_3 d_q & -x_1 d_q + x_2 + x_3 c_q \\
-x_1 d_q - x_2 + x_3 c_q  & x_0 - x_1 c_q - x_3 d_q \\
\end{pmatrix}\right],
\]
where $x = x_0 + x_1 i + x_2 j + x_3 k \in \mathbb{H}(\mathbb{Q})^{\ast}$.
The following homomorphism $\psi_{p,l}$
will play a crucial role in our analysis of quaternion groups.
Let 
\[
\psi_{p,l} : \mathbb{H}(\mathbb{Q})^{\ast} \to 
\mathrm{PGL}_2(\mathbb{Q}_p) \times \mathrm{PGL}_2(\mathbb{Q}_l)
\]
be given by 
\[
\psi_{p,l}(x) := (\psi_p(x), \psi_l(x)).
\]
This homomorphism is not injective, 
in fact (see \cite[Chapter~3]{Rattaggi})
\[
\mathrm{ker}(\psi_{p,l}) = Z(\mathbb{H}(\mathbb{Q})^{\ast}) =
\{x \in \mathbb{H}(\mathbb{Q})^{\ast} : x = \overline{x}\} \cong \mathbb{Q}^{\ast},
\]
where for the last isomorphism we identify $\mathbb{Q}^{\ast}$
with the image of the natural injective homomorphism 
$\mathbb{Q}^{\ast} \to \mathbb{H}(\mathbb{Q})^{\ast}$
given by $x_0 \mapsto x_0 + 0 \cdot i + 0 \cdot j + 0 \cdot k$.
In particular, we have $\psi_{p,l}(x) = \psi_{p,l}(y)$, if and only if
$y = \lambda x$ for some $\lambda \in \mathbb{Q}^{\ast}$, and
therefore
\[
\psi_{p,l}(x) = \psi_{p,l}(-x)
\] 
for each $x \in \mathbb{H}(\mathbb{Q})^{\ast}$.
Moreover, using the rule $\overline{x} = |x|^2 x^{-1}$, we also get 
\[
\psi_{p,l}(\overline{x}) = \psi_{p,l}(x)^{-1}.
\]

For an odd prime number $q$, let $X_q$ be the finite set of integral quaternions
\begin{align*}
X_q := \{x = x_0  + &x_1 i + x_2 j + x_3 k \in \mathbb{H}(\mathbb{Z}) \,; 
\quad |x|^2 = q \, ; \notag \\
&x_0 \text{ odd}, x_1, x_2, x_3 \text{ even}, \text{ if } q \equiv 1
\!\!\!\! \pmod 4\,; \notag \\
&x_1 \text{ even}, x_0, x_2, x_3 \text{ odd}, \text{ if } q \equiv 3
\!\!\!\! \pmod 4
\} \, .
\end{align*}
Observe that $X_q$ has exactly $2(q+1)$ elements 
(by Jacobi's theorem on the number of representations of an integer as a sum of four squares)
and that $X_q$ is closed under conjugation and under multiplication by $-1$.
As examples we have 
\[
X_3 = \{ \pm 1 \pm j \pm k \},
\] 
where all of the $2^3$ possible combinations of signs are allowed,
and 
\[
X_5 = \{ \pm 1 \pm 2i, \, \pm 1 \pm 2j, \, \pm 1 \pm 2k \}.
\]

Finally, let $Q_{p,l}$ be the subgroup of $\mathbb{H}(\mathbb{Q})^{\ast}$
generated by $(X_p \cup X_l) \subset \mathbb{H}(\mathbb{Z})$ and let 
$\Gamma_{p,l} < \mathrm{PGL}_2(\mathbb{Q}_p) \times \mathrm{PGL}_2(\mathbb{Q}_l)$ 
be its image $\psi_{p,l}(Q_{p,l})$.
Using the properties $\psi_{p,l}(x) = \psi_{p,l}(-x)$ and
$\psi_{p,l}(\overline{x}) = \psi_{p,l}(x)^{-1}$ mentioned above,
it follows that $\Gamma_{p,l}$ is generated by $(p+1)/2 + (l+1)/2$
elements.
For example $\Gamma_{3,5}$ is generated by
the five elements $\psi_{3,5}(1+j+k)$, $\psi_{3,5}(1+j-k)$,
$\psi_{3,5}(1+2i)$, $\psi_{3,5}(1+2j)$ and $\psi_{3,5}(1+2k)$.
See Section~\ref{SectionEx} for a finite presentation of $\Gamma_{3,5}$.

We recall some known properties of the group $\Gamma_{p,l}$
from \cite{BMII}, \cite{Mozes1}, \cite{Rattaggi}, \cite{Rattaggi2}, \cite{RaRo}.
It has a finite presentation with generators 
\[
a_1, \ldots, a_{\frac{p+1}{2}}, b_1, \ldots, b_{\frac{l+1}{2}}
\]
and $(p+1)(l+1)/4$ defining relations of the form
$ab\tilde{a}\tilde{b} = 1$ for some
\[
a, \tilde{a} \in \{a_1, \ldots, a_{\frac{p+1}{2}}\}^{\pm 1} 
\text{\, and \,}
b, \tilde{b} \in \{b_1, \ldots, b_{\frac{l+1}{2}}\}^{\pm 1}.
\]
It acts freely and transitively on the vertices of the product
of two regular trees of degree $p+1$ and $l+1$, respectively, is CAT(0), bi-automatic
and can be decomposed as an amalgamated product of finitely generated free groups
(see \cite{BMII} or \cite{Rattaggi}).
It is CSA, i.e.\ all maximal abelian subgroups are malnormal (see \cite[Proposition~2.6]{RaRo}),
in particular it is commutative transitive, i.e.\ the relation of commutativity is 
transitive on non-trivial elements.
It is also linear (see \cite[Proposition~31]{Rattaggi2} for an explicit injective homomorphism
$\Gamma_{p,l} \to \mathrm{SO}_3(\mathbb{Q}) < \mathrm{GL}_3(\mathbb{Q}))$,
contains free abelian subgroups $\mathbb{Z} \times \mathbb{Z}$ (see \cite[Proposition~4.2(3)]{Rattaggi})
as well as non-abelian free subgroups, for example
\[
\langle a_1, \ldots, a_{\frac{p+1}{2}} \rangle \cong F_{\frac{p+1}{2}}
\text{ or } 
\langle b_1, \ldots, b_{\frac{l+1}{2}} \rangle \cong
F_{\frac{l+1}{2}}.
\]
(However, it is not known whether there are elements
\[
a \in \langle a_1, \ldots, a_{\frac{p+1}{2}} \rangle \text{\, and \,} 
b \in \langle b_1, \ldots, b_{\frac{l+1}{2}} \rangle
\]
generating a free group $\langle a,b \rangle \cong F_2$.)
The conditions on the parity of $x_0$, $x_1$, $x_2$, $x_3$ in the definition of
$X_p$ and $X_l$ are mainly used to guarantee that $\Gamma_{p,l}$ is torsion-free
(see \cite[Proposition~3.6]{Mozes1} and
\cite[Theorem~3.30(4)]{Rattaggi}).
Since finitely generated, torsion-free,
virtually free groups are free
(see \cite{Stallings}),
the group $\Gamma_{p,l}$ is not virtually free.
Moreover,
$\Gamma_{p,l}$ is not virtually abelian. 
(It is well-known that the property of being virtually abelian is invariant under 
quasi-isometry for finitely generated groups.
Now we use that $\Gamma_{p,l}$ is quasi-isometric to the non-virtually abelian group $F_2 \times F_2$, 
see \cite[Proposition~4.25(4)]{Rattaggi}. 
Alternatively, without using finite index subgroups of 
$F_2 \times F_2$, we note that $\Gamma_{p,l}$ is also
quasi-isometric to the non-abelian finitely presented torsion-free simple groups described
in \cite{BMII} and \cite{Rattaggi}. Infinite non-abelian simple groups
are obviously not virtually abelian,
which gives another proof that $\Gamma_{p,l}$ is not virtually abelian.) 
Any non-trivial normal subgroup of $\Gamma_{p,l}$ has finite index
(by the ``Normal Subgroup Theorem'' of Burger-Mozes
\cite[Chapter~4~and~5]{BMII}). 
This property also holds for any finite index subgroup of $\Gamma_{p,l}$.

\section{The non-commutative case}
We begin with some basic notations and a general lemma which will be applied
to our quaternion groups later.
If $G$ is any group and $g_1, g_2 \in G$ two elements,
we denote by $[g_1,g_2] := g_1 g_2 g_1^{-1} g_2^{-1}$
the commutator of $g_1$ and $g_2$, by $G'$ the commutator subgroup of $G$,
by $G^{ab}$ the quotient $G/G'$, and by $ZG$ or $Z(G)$ the center of $G$.

\begin{Lemma} \label{Lemma1}
Let $G$ be a (multiplicatively written) group and $N$ a normal subgroup of $G$.
Then we have the following commutative diagram with exact
rows and exact columns
\[
\xymatrix{
          &    1                         &  1                          &  1                         &   \\
1 \ar[r]  & (G/N)' \ar[r]^{i_1} \ar[u]   & G/N \ar[r]^{p_1} \ar[u]     & (G/N)^{ab} \ar[r] \ar[u]   & 1 \\
1 \ar[r]  & G' \ar[r]^{i_2} \ar[u]^{q_1} & G \ar[r]^{p_2} \ar[u]^{q_2} & G^{ab} \ar[r] \ar[u]^{q_3} & 1 \\
1 \ar[r]  & N \cap G' \ar[r]^{i_3} \ar[u]^{j_1} & N \ar[r]^{p_3} \ar[u]^{j_2} & p_2(N) \ar[r] \ar[u]^{j_3} & 1 \\
          & 1 \ar[u]                     &  1 \ar[u]                   &  1 \ar[u]                  &
} 
\]
where $i_1$, $i_2$, $i_3$, $j_1$, $j_2$, $j_3$ are the natural homomorphisms injecting the corresponding normal subgroups,
$p_1$, $p_2$, $q_2$ are the natural projections, $p_3$ is the restriction of $p_2$ to $N$,
$q_1$ is induced by 
\[
q_1([g_1,g_2]) := [g_1N,g_2N] = [g_1,g_2]N,
\] 
i.e.\
$q_1(g') = g'N$, if $g' \in G'$, and finally $q_3(gG') := gN(G/N)'$.
\end{Lemma}

\begin{proof}
It is clear that the top row, the middle row and the middle column are exact.
The exactness of the bottom row and left column also follows immediately, since
\[
\mathrm{ker}(p_3) = \{n \in N : nG' = G'\} = i_3(N \cap G')
\]
and
\[
\mathrm{ker}(q_1) = \{g' \in G' : g'N = N\} = j_1(N \cap G').
\]
It remains to see that the right column is exact.
The only non-obvious part is to show that $\mathrm{ker}(q_3) = \mathrm{im}(j_3)$.
We have by definition of $q_3$
\[
\mathrm{ker}(q_3) = \{ gG' \in G/G' : gN \in (G/N)' \}.
\]
Since
\begin{align}
gN \in (G/N)' \, &\Leftrightarrow \, gN \in \{g'N \in G/N : g' \in G' \} \,
               \Leftrightarrow \, \exists \, g' \in G' : gN = g'N \notag \\
              &\Leftrightarrow \exists \, g' \in G' : g'^{-1}g \in N \,
               \Leftrightarrow \, \exists \, g' \in G', n \in N  : g'^{-1}g = n \notag \\
              &\Leftrightarrow \, \exists \, g' \in G', n \in N  : g'^{-1} = ng^{-1} \,
               \Leftrightarrow \, \exists \, n \in N  : ng^{-1} \in G' \notag \\
              &\Leftrightarrow \, \exists \, n \in N  : g^{-1}G' = n^{-1}G' \,
               \Leftrightarrow \, \exists \, n \in N  : gG' =  nG' \notag \\
              &\Leftrightarrow \, gG' \in \{ nG' : n \in N \} \notag
\end{align}
it follows that 
\[
\mathrm{ker}(q_3) = \{ nG' \in G/G' : n \in N \} = \mathrm{im}(j_3).
\]
The commutativity of the diagram is a direct consequence of the given definitions.
\end{proof}

We refer to Section~\ref{SectionPrel} for the notations
concerning our quaternion groups.
Let $x, y \in Q_{p,l}$ be two non-commuting quaternions and 
$\psi$ the restriction of $\psi_{p,l}$ 
to the subgroup $\langle x, y \rangle < Q_{p,l} < \mathbb{H}(\mathbb{Q})^{\ast}$.
We want to apply Lemma~\ref{Lemma1} to the situation where
$q_2 : G \to G/N$ is the surjective homomorphism 
$\psi : \langle x, y \rangle \to \langle \psi(x), \psi(y) \rangle$ and where
$N = \mathrm{ker}(\psi)$, in order to
get some information on the structure of the four
groups $\langle x, y \rangle$, $Z\langle x, y \rangle$, $\langle x, y \rangle'$ and $\langle x, y \rangle^{ab}$.
First, we investigate the group $\mathrm{ker}(\psi)$ in
Lemma~\ref{Lemma3}(1),(2),(3), applying the following lemma.
Then, we try to understand in Lemma~\ref{Lemma3}(4) the bottom left term $N \cap G'$
of the diagram.

\begin{Lemma} \label{Lemma2}
\begin{itemize}
\item[(1)]
The group $\mathbb{H}(\mathbb{Q})^{\ast}$ is commutative transitive
on non-central elements, in other words
$xz = zx$, $yz = zy$ implies $xy = yx$,
whenever $x, y, z \in \mathbb{H}(\mathbb{Q})^{\ast} \setminus
\mathbb{Q}^{\ast}$.
\item[(2)]
The group $\mathbb{H}(\mathbb{Q})^{\ast}$ 
contains no subgroup isomorphic to a direct product of two 
non-abelian groups.
\end{itemize}
\end{Lemma}

\begin{proof}
\begin{itemize}
\item[(1)]
This is an elementary computation, see \cite[Lemma~3.4(3)]{Rattaggi}.
\item[(2)]
Suppose that $\mathbb{H}(\mathbb{Q})^{\ast}$ contains a subgroup 
$G \times H$, where $G, H$ are non-abelian groups.
Take $g_1, g_2 \in G$ such that $g_1 g_2 \ne g_2 g_1$
and $h_1, h_2 \in H$ such that $h_1 h_2 \ne h_2 h_1$.
Then clearly $g_1, g_2, h_1 \notin Z(\mathbb{H}(\mathbb{Q})^{\ast})$.
The fact that $h_1$ commutes with $g_1$ and with $g_2$ but $g_1 g_2 \ne g_2 g_1$
now contradicts part~(1) of this lemma.
\end{itemize}
\end{proof}

\begin{Lemma} \label{Lemma3}
Let $x, y \in Q_{p,l}$ be two non-commuting quaternions and
$\psi$ the restriction of $\psi_{p,l}$ 
to $\langle x, y \rangle$. Then
\begin{itemize}
\item[(1)] $Z \langle x, y \rangle = \langle x, y \rangle \cap \mathbb{Q}^{\ast}$.
\item[(2)] $\mathrm{ker}(\psi) = Z \langle x, y \rangle$.
\item[(3)] $Z \langle x, y \rangle < \{ \pm p^r l^s : r,s \in \mathbb{Z}\} 
= \langle -1, p, l \rangle
\cong \mathbb{Z}_2 \times \mathbb{Z} \times \mathbb{Z}$,
in particular $Z \langle x, y \rangle$ is finitely presented.
\item[(4)] $Z \langle x, y \rangle \cap \langle x, y \rangle' < \{\pm 1\} \cong \mathbb{Z}_2$.
\item[(5)] The following diagram is commutative and has exact rows and columns
\[
\xymatrix
{
          &    1                         &  1                          &  1                         &   \\
1 \ar[r]  & \langle \psi(x), \psi(y) \rangle' \ar[r] \ar[u]   & \langle \psi(x), \psi(y) \rangle 
\ar[r] \ar[u]     & \langle \psi(x), \psi(y) \rangle^{ab} \ar[r] \ar[u]   & 1 \\
1 \ar[r]  &  \langle x, y \rangle' \ar[r] \ar[u] & \langle x, y \rangle \ar[r]^{p_2} 
\ar[u]^{\psi} &\langle x, y \rangle^{ab} \ar[r] \ar[u] & 1 \\
1 \ar[r]  & 1 \text{ or } \{\pm 1\} \ar[r] \ar[u] & Z\langle x, y \rangle 
\ar[r] \ar[u] & p_2(Z\langle x, y \rangle) \ar[r] \ar[u] & 1 \\
          & 1 \ar[u]                     &  1 \ar[u]                   &  1 \ar[u]                  &
} 
\]
where the maps are defined as in Lemma~\ref{Lemma1},
starting with
\[
(N \longrightarrow G \overset{q_2}{\longrightarrow} G/N) \, := \,
(\mathrm{ker}(\psi) \longrightarrow \langle x, y \rangle \overset{\psi}{\longrightarrow} 
\langle \psi(x), \psi(y) \rangle).
\]
\item[(6)] $[x,y] \ne -1$.
\item[(7)]  Define the following three statements:\\
(i) $Z\langle x, y \rangle$ is infinite. \\
(ii) $|x|^2 \ne 1$ or $|y|^2 \ne 1$. \\
(iii) $\langle x, y \rangle^{ab}$ is infinite.\\
Then (i) implies (ii) which in turn implies (iii).
\item[(8)] Let $r(\psi(x),\psi(y))$ be a relator in the
group generated by $\psi(x)$, $\psi(y)$.
Then $r(x,y) \in Z\langle x, y \rangle$.
\item[(9)] $\langle x,y \rangle \cong F_2$ if and only if 
$\langle \psi(x),\psi(y) \rangle \cong F_2$.
\end{itemize}
\end{Lemma}

\begin{proof}
\begin{itemize}
\item[(1)] First recall that $Z(\mathbb{H}(\mathbb{Q})^{\ast}) = \mathbb{Q}^{\ast}$.
Since $x$, $y$ do not commute by assumption, we have in particular $x, y \notin \mathbb{Q}^{\ast}$.
If $z \in Z \langle x, y \rangle$, then $z \in \mathbb{Q}^{\ast}$
(otherwise $x$, $y$, $z$ would be three non-central elements in $\mathbb{H}(\mathbb{Q})^{\ast}$
such that $z$ commutes with $x$ and $y$, but $x$ and $y$ do not
commute, contradicting Lemma~\ref{Lemma2}(1)).
This shows that 
$Z \langle x, y \rangle < \langle x, y \rangle \cap \mathbb{Q}^{\ast}$.

To show $Z \langle x, y \rangle > \langle x, y \rangle \cap \mathbb{Q}^{\ast}$
we again use that $\mathbb{Q}^{\ast} = Z(\mathbb{H}(\mathbb{Q})^{\ast})$.
\item[(2)] Using $\mathrm{ker}(\psi_{p,l}) = \mathbb{Q}^{\ast}$ 
and part~(1) of this lemma, we get
\[
\mathrm{ker}(\psi) = \mathrm{ker}(\psi_{p,l}) \cap \langle x, y \rangle = 
\mathbb{Q}^{\ast} \cap \langle x, y \rangle = Z \langle x, y \rangle.
\]
\item[(3)] The norm of any element in $\langle x, y \rangle$
is of the form $p^r l^s$ for some $r,s \in \mathbb{Z}$, hence
\[
Z \langle x, y \rangle = \langle x, y \rangle \cap \mathbb{Q}^{\ast}
= \langle x, y \rangle \cap \{ \pm p^r l^s : r,s \in \mathbb{Z}\}
< \{ \pm p^r l^s : r,s \in \mathbb{Z}\}.
\]
The center $Z \langle x, y \rangle$ is finitely presented, since 
subgroups of finitely generated abelian groups
are again finitely generated abelian, 
and since finitely generated abelian groups are finitely presented.
\item[(4)] As seen in part~(1) of this lemma, we have
$Z \langle x, y \rangle < \mathbb{Q}^{\ast}$.
By the multiplicativity of the norm, any commutator in $\mathbb{H}(\mathbb{Q})^{\ast}$
has norm $1$, hence any element in $\langle x, y \rangle'$ has norm~$1$.
The only elements in $\mathbb{Q}^{\ast}$ having norm $1$ are $\pm 1$,
and the statement follows.
\item[(5)] We combine Lemma~\ref{Lemma1}
with part~(2) and (4) of this lemma.
\item[(6)] $[x,y] = -1$ is equivalent to $xy = -yx$ which implies $\Re(x) = 0$ 
(see \cite[Lemma~3.4(2)]{Rattaggi}).
But then 
\[
x^2 = -x_1^2 -x_2^2 -x_3^2 \in \mathbb{Q}^{\ast}, 
\]
hence
\[
1_{\Gamma_{p,l}} = \psi(x^2) = \psi(x)^2.
\]
Since $\Gamma_{p,l}$ is torsion-free, $\psi(x) = 1_{\Gamma_{p,l}}$,
and $x \in \mathrm{ker}(\psi) = \mathbb{Q}^{\ast} \cap \langle x, y \rangle$
by part~(1) and (2) of this lemma,
contradicting $\Re(x) = 0$.
\item[(7)]  
If $|x|^2 = |y|^2 = 1$, then all elements of $\langle x, y \rangle$ have norm $1$,
in particular all elements of $Z\langle x, y \rangle$ have norm $1$.
As in the proof of part~(4) of this lemma we see
that $Z\langle x, y \rangle < \{\pm 1\}$ is finite.

To show that (ii) implies (iii),
we assume without loss of generality that $|x|^2 \ne 1$.
Then the norms of $x^m$ and $x^n$ are distinct, whenever $m \ne n \in \mathbb{Z}$.
Any element in $\langle x, y \rangle'$ has norm $1$, as observed in the
proof of part~(4) of this lemma.
It follows that any two cosets $x^m \langle x, y \rangle'$ and $x^n \langle x, y \rangle'$
are distinct.
In particular, $\langle x, y \rangle'$ has infinite index in $\langle x, y \rangle$,
hence $\langle x, y \rangle^{ab}$ is infinite.
\item[(8)] We have $r(x,y) \in \mathrm{ker}(\psi)$ and apply part~(2) of this lemma.
\item[(9)] See \cite[Proposition~32]{Rattaggi2}.
\end{itemize}
\end{proof}

\begin{Remark} \label{Rem3b}
Going through the proofs of Lemma~\ref{Lemma3} we see that
some results can be easily generalized:
Let $G$ be any non-abelian subgroup of
$\mathbb{H}(\mathbb{Q})^{\ast}$.
Then
\begin{itemize}
\item $\mathrm{ker}(\psi_{p,l}|_G) = ZG = G \cap \mathbb{Q}^{\ast}$.
\item $ZG \cap G' < \{\pm 1\} $.
\item $ZG$ is infinite $\Rightarrow$ $G \nless \mathbb{H}(\mathbb{Q})_1$
$\Rightarrow$ $G^{ab}$ is infinite.
\end{itemize}
\end{Remark}

Related to Lemma~\ref{Lemma3}(4),
note that clearly $Z \langle x, y \rangle \cap \langle x, y \rangle' = 1$
if $x, y \in Q_{p,l}$ commute.
We conjecture that 
$-1 \in Z \langle x, y \rangle \cap \langle x, y \rangle'$ 
is also impossible if $x, y \in Q_{p,l}$ do not commute (see Conjecture~\ref{Conj6}).

\begin{Question} \label{Question4}
Let $x, y \in Q_{p,l}$ be two non-commuting quaternions.
Is it possible that $-1 \in Z \langle x, y \rangle$?
(Equivalently, is it possible that $-1 \in \langle x, y \rangle$?)
\end{Question}

\begin{Conjecture} \label{Conj5}
Let $x, y \in Q_{p,l}$.
Then $-1 \notin \langle x, y \rangle'$.
We even conjecture that $-1 \notin Q_{p,l}'$ (cf.\ Remark~\ref{Rem14}).
\end{Conjecture}

As a consequence of Lemma~\ref{Lemma3}(4) and Conjecture~\ref{Conj5} we have
\begin{Conjecture} \label{Conj6}
Let $x, y \in Q_{p,l}$. Then
$Z \langle x, y \rangle \cap \langle x, y \rangle' = 1$.
\end{Conjecture}

In contrast to Lemma~\ref{Lemma3}(6), commutators in
$\mathbb{H}(\mathbb{Q})^{\ast}$ can of course be $-1$, for example
$[i,j] = -1$. 
More generally, if $x, y \in \mathbb{H}(\mathbb{Q})^{\ast}$, 
then $[x,y] = -1$ if and only if $x_0 = y_0 = 0$
and $x_1y_1 + x_2 y_2 + x_3y_3 = 0$, see \cite[Lemma~1]{Rattaggi3}. 
In particular we have
$-1 \in ({\mathbb{H}(\mathbb{Q})^{\ast}})'$
and $Z ({\mathbb{H}(\mathbb{Q})^{\ast}}) \cap ({\mathbb{H}(\mathbb{Q})^{\ast}})' = \{\pm 1\}$.
However, it is easy to proof that $-1 \notin \langle X_q \rangle'$
for all odd prime numbers $q$.

\begin{Problem} \label{Problem7}
Characterize the (non-abelian) subgroups $G$ of $\mathbb{H}(\mathbb{Q})^{\ast}$
such that $-1 \in G'$.
\end{Problem}

Before applying Lemma~\ref{Lemma3} to deduce our main results of
this section in Theorem~\ref{Thm11},
we prove another general lemma and
some statements about groups generated by quaternions $x, y \in Q_{p,l}$ 
satisfying a (relatively weak) norm condition in Proposition~\ref{Prop9}.

\begin{Lemma} \label{Lemma8}
The only non-trivial element of finite order in $Q_{p,l}$ is $-1$.
In particular, a subgroup $G < Q_{p,l}$ is torsion-free if and
only if $-1 \notin G$.
\end{Lemma}

\begin{proof}
Let $z \in Q_{p,l}$ such that $z^n = 1$ for some $n \in \mathbb{N}$.
Then $\psi(z)^n = 1_{\Gamma_{p,l}}$.
Since $\Gamma_{p,l}$ is torsion-free, we have $\psi(z) = 1_{\Gamma_{p,l}}$,
hence $z \in \mathbb{Q}^{\ast}$.
But then $z \in \{ \pm 1 \}$, since 
\[
1 = |1|^2 = |z^n|^2 = (|z|^2)^n = z^{2n}.
\]
\end{proof}

\begin{Proposition} \label{Prop9}
Let $x, y \in Q_{p,l}$ be two quaternions
of norms $|x|^2 = p^{r_1}l^{s_1}$,
$|y|^2 = p^{r_2}l^{s_2}$, $r_1, r_2, l_1, l_2 \in \mathbb{Z}$,
such that $r_1 s_2 \ne r_2 s_1$.
(This condition holds for example if $|x|^2 = p^r$, $|y|^2 = l^s$ for some
$r, s \in \mathbb{Z} \setminus \{ 0 \}$).
Then
\begin{itemize}
\item[(1)] $\langle x, y \rangle^{ab} \cong \mathbb{Z} \times \mathbb{Z}$,
generated by the two commuting elements 
$x \langle x,y \rangle'$ and $y \langle x,y \rangle'$.
\item[(2)] $\langle x, y \rangle \cap \mathbb{H}(\mathbb{Q})_1 = \langle x, y \rangle'$.
\item[(3)] $\langle x, y \rangle$ is torsion-free if and only if 
$\langle x, y \rangle'$ is torsion-free.
\end{itemize}
\end{Proposition}

\begin{proof}
\begin{itemize}
\item[(1)] Let $r(x,y) = 1$ be any relation in $\langle x, y \rangle$.
Let $s_x$ be the exponent sum of $x$ in $r(x,y)$
and $s_y$ the exponent sum of $y$ in $r(x,y)$.
Taking the norm of $r(x,y) = 1$, we get
\[
p^{s_x r_1 + s_y r_2} \cdot l^{s_x s_1 + s_y s_2} = 1,
\]
hence
\[
\left( \begin{array}{cc}
r_1  &  r_2  \\
s_1  &  s_2  \notag
\end{array} \right)
\begin{pmatrix}
s_x \\
s_y 
\end{pmatrix}
=
\begin{pmatrix}
0 \\
0 
\end{pmatrix}.
\]
By assumption, the determinant $r_1 s_2 - r_2 s_1$ is
non-zero and we get $s_x = s_y = 0$.
This shows that the relator $r(x,y)$ is a consequence of
the relator $[x,y]$.
Thus $\langle x, y \rangle^{ab} \cong \langle x,y \mid [x,y] \rangle$,
where the two generators correspond to 
$x \langle x,y \rangle'$ and $y \langle x,y \rangle'$.
(Note that for any group $G$, a presentation of $G^{ab}$
is obtained from a presentation $\langle X \mid R \rangle$ of $G$ 
by adding to $R$ all 
commutators of elements of $X$.
In the case of a $2$-generator group, we therefore have to 
add just a single commutator.) 
\item[(2)] Clearly $\langle x, y \rangle \cap \mathbb{H}(\mathbb{Q})_1  > \langle x, y \rangle'$.

To show the other inclusion, let $z \in \langle x, y \rangle \cap \mathbb{H}(\mathbb{Q})_1$.
Since $z$ has norm~$1$, it follows as in part~(1) of this proposition that 
the exponent sums of $x$ and $y$ in $z$ are $0$.
Thus $z \langle x, y \rangle' = \langle x, y \rangle'$
in the abelian group $\langle x, y \rangle / \langle x, y \rangle'$,
in other words $z \in \langle x, y \rangle'$.
\item[(3)]
Suppose that $\langle x, y \rangle$ is not torsion-free.
Then $-1 \in \langle x, y \rangle$ by Lemma~\ref{Lemma8},
hence $-1 \in \langle x, y \rangle'$ by part~(2) of this
proposition and $\langle x, y \rangle'$ is not torsion-free. 

The other direction is obvious.
\end{itemize}
\end{proof}

\begin{Lemma} \label{Lemma10}
Let $x, y \in Q_{p,l}$ be two quaternions and assume that
$\langle \psi(x), \psi(y) \rangle$ has finite index in $\Gamma_{p,l}$.
Then
\begin{itemize}
\item[(1)] $x$ and $y$ do not commute.
\item[(2)] $\langle \psi(x), \psi(y) \rangle^{ab}$ is finite.
\end{itemize}
\end{Lemma}

\begin{proof}
\begin{itemize}
\item[(1)] If $x$ and $y$ commute, then
$\langle \psi(x), \psi(y) \rangle$ is an abelian group,
but $\Gamma_{p,l}$ is not virtually abelian.
\item[(2)] If $G$ is a subgroup of finite index in $\Gamma_{p,l}$,
then any non-trivial normal subgroup of $G$ has finite index in $G$.
Since $\langle \psi(x), \psi(y) \rangle$ is
not abelian, the normal subgroup $\langle \psi(x), \psi(y) \rangle'$
is not trivial, hence has finite index in $\langle \psi(x), \psi(y) \rangle$.
\end{itemize}
\end{proof}

\begin{Theorem} \label{Thm11}
Let $x, y \in Q_{p,l}$ be two quaternions and let $\psi$ be the
restriction of $\psi_{p,l}$ to $\langle x, y \rangle$. Assume that the group
$\psi(\langle x, y \rangle) = \langle \psi(x), \psi(y) \rangle$ has finite index in $\Gamma_{p,l}$.
Then
\begin{itemize}
\item[(1)] The group $\langle x, y \rangle$ is finitely presented and infinite.
\item[(2)] The group $\langle x, y \rangle'$ is finitely presented and infinite.
\item[(3)] The following three statements are equivalent: \\
(i) $Z\langle x, y \rangle$ is infinite. \\
(ii) $|x|^2 \ne 1$ or $|y|^2 \ne 1$. \\
(iii) $\langle x, y \rangle^{ab}$ is infinite.
\item[(4)] The group $\langle x, y \rangle$ is not virtually solvable.
\item[(5)] The group $\langle x, y \rangle$ contains a free subgroup $F_2$ of infinite index.
\item[(6)] The group $\langle x, y \rangle$ contains a subgroup $\mathbb{Z} \times \mathbb{Z}$
of infinite index. In particular, $\langle x, y \rangle$ is not hyperbolic.
\item[(7)] Let $N$ be a normal subgroup of $\langle x, y \rangle$ of infinite index
such that $Z\langle x, y \rangle < N$. Then $N = Z\langle x, y \rangle$.
\end{itemize}
\end{Theorem}

\begin{proof}
By Lemma~\ref{Lemma10}(1) we can use the commutative diagram
of Lemma~\ref{Lemma3}(5)
which we will simply call ``the diagram'' in this proof.
\begin{itemize}
\item[(1)] The group $\langle \psi(x), \psi(y) \rangle$ is finitely presented,
since it has finite index in the finitely presented group $\Gamma_{p,l}$
by assumption.
Using the middle column of the diagram,
the group $\langle x, y \rangle$ is an extension of the finitely presented group
$Z\langle x, y \rangle$ by the finitely presented group $\langle \psi(x), \psi(y) \rangle$,
hence finitely presented.

It is clear that $\langle x, y \rangle$ is infinite, since $\langle \psi(x), \psi(y) \rangle$ is infinite
as a finite index subgroup of the infinite group $\Gamma_{p,l}$.
\item[(2)] By Lemma~\ref{Lemma10}(2) $\langle \psi(x), \psi(y) \rangle^{ab}$ is finite.
By the exactness of the top row in the diagram,
$\langle \psi(x), \psi(y) \rangle'$ has finite index in $\langle \psi(x), \psi(y) \rangle$,
hence is also finitely presented.
Now, using the exactness of the left column in the diagram,
$\langle x, y \rangle'$ is an extension of a finite group ($\mathbb{Z}_2$ or $1$) by
the finitely presented group $\langle \psi(x), \psi(y) \rangle'$ and therefore finitely presented.

By the exactness of the top row and left column,
$\langle \psi(x), \psi(y) \rangle'$ and $\langle x, y \rangle'$ are
infinite groups.
\item[(3)] Because of Lemma~\ref{Lemma10}(1) and Lemma~\ref{Lemma3}(7),
it remains to prove that (iii) implies (i).
Therefore suppose that $\langle x, y \rangle^{ab}$ is infinite.
Since $\langle \psi(x), \psi(y) \rangle^{ab}$ is finite by
Lemma~\ref{Lemma10}(2),
the exactness of the right column in the diagram shows
that $p_2(Z\langle x, y \rangle)$ is infinite.
The exactness of the bottom row in the diagram
now shows that $Z\langle x, y \rangle$ is infinite.
\item[(4)] We first show that the group $\langle \psi(x), \psi(y) \rangle$ is not virtually solvable.
Note that the property of being virtually solvable is \emph{not} invariant
under quasi-isometry for finitely generated groups (see \cite{Dyubina}),
why we cannot use the same strategy as for the proof that $\Gamma_{p,l}$ is not virtually abelian.
Instead of that, we adapt an idea already used in the proof of \cite[Corollary~2.8]{RaRo}.
Let $V < \langle \psi(x), \psi(y) \rangle$ be a subgroup of finite index.
Then $V$ is not abelian (otherwise $\Gamma_{p,l}$ would be virtually abelian).
The group $\langle \psi(x), \psi(y) \rangle$ is a non-abelian CSA-group,
since $\Gamma_{p,l}$ is CSA, and subgroups of CSA-groups are CSA (\cite[Proposition~10(3)]{MR}).
Now, we use the fact that non-abelian CSA-groups do not have any non-abelian solvable
subgroups (\cite[Proposition~9(5)]{MR}) to conclude that $V$ is not solvable.

Let $U$ be any finite index subgroup of $\langle x,y \rangle$.
Since quotients of solvable groups are solvable, and
\[
[\langle \psi(x), \psi(y) \rangle : \psi(U)] \leq [\langle x,y \rangle : U] < \infty,
\]
the group $U$ is not solvable, and therefore $\langle x,y \rangle$ is not
virtually solvable.
\item[(5)]
There is a well-known injective homomorphism of groups
\begin{align}
\mathbb{H}(\mathbb{Q})^{\ast} &\to \mathrm{GL}_4(\mathbb{Q}) \notag \\
x_0 + x_1 i + x_2 j + x_3 k   &\mapsto
\left(
\begin{array}{rrrr}
x_0 & -x_1 & -x_2 & -x_3 \\
x_1 &  x_0 & -x_3 &  x_2 \\
x_2 &  x_3 &  x_0 & -x_1 \\
x_3 & -x_2 &  x_1 &  x_0 \\
\end{array}
\right), \notag
\end{align}
in particular $\langle x,y \rangle$ is a finitely generated linear group in characteristic $0$.
Since $\langle x,y \rangle$ is not virtually solvable by part~(4) of this theorem, 
it contains by the Tits Alternative (\cite{Tits}) a free subgroup $F_2$.

A free subgroup $F_2$ of $\langle x,y \rangle$ always has infinite index, since otherwise
$\psi(F_2)$ has finite index in $\Gamma_{p,l}$, but
$\psi(F_2) \cong F_2$ by Lemma~\ref{Lemma3}(9) and
$\Gamma_{p,l}$ is not virtually free.
\item[(6)] Let $a_1, \ldots, a_{\frac{p+1}{2}}, b_1, \ldots, b_{\frac{l+1}{2}}$ be the standard
generators of $\Gamma_{p,l}$.
By \cite[Proposition~4.2(3)]{Rattaggi}, there are elements 
\[
a \in \langle a_1, \ldots, a_{\frac{p+1}{2}} \rangle
\text{\, and \,} 
b \in \langle b_1, \ldots, b_{\frac{l+1}{2}} \rangle,
\]
such that $a$, $b$ generate a subgroup of $\Gamma_{p,l}$ isomorphic to $\mathbb{Z} \times \mathbb{Z}$.
Since $\langle \psi(x), \psi(y) \rangle$ has finite index in
$\Gamma_{p,l}$ by assumption,
there are $m, n \in \mathbb{N}$ such that $a^m, b^n \in \langle \psi(x), \psi(y) \rangle$
and 
\[
\mathbb{Z} \times \mathbb{Z} \cong \langle a^m, b^n \rangle < \langle \psi(x), \psi(y) \rangle.
\]
Let $w_1$, $w_2$ be quaternions in $\langle x, y \rangle$ 
such that $a^m = \psi(w_1)$ and $b^n = \psi(w_2)$.
Since $a^m$ commutes with $b^n$,
we have 
\[
\psi(w_1 w_2)  = \psi(w_2 w_1),
\]  
hence either 
\[
w_1 w_2 = -w_2 w_1
\]
or 
\[
w_1 w_2 = w_2 w_1.
\]
The first case can be excluded similarly as in Lemma~\ref{Lemma3}(6).
The commuting quaternions $w_1$, $w_2$ have infinite order,
since $w_1, w_2 \notin \{ -1, 1 \}$.
Moreover, it is not possible that
$w_1^r = w_2^s$ for some $r, s \in \mathbb{N}$,
since otherwise $\psi(w_1)^r = \psi(w_2)^s$, i.e.\ $a^{mr} = b^{ns}$.
But since $mr \ne 0 \ne ns$, 
this is impossible using the ``normal form'' (see \cite[Proposition~1.10]{Rattaggi})
for elements of $\Gamma_{p,l}$.
Therefore 
\[
\mathbb{Z} \times \mathbb{Z} \cong \langle w_1, w_2 \rangle < \langle x, y \rangle.
\]
This subgroup has to be of infinite index, since $\Gamma_{p,l}$ is not virtually abelian.
\item[(7)] The assumptions imply that $N/Z\langle x,y \rangle$ is a normal subgroup
of $\langle x,y \rangle / Z\langle x,y \rangle$ of infinite index,
since 
\[
(\langle x,y \rangle / Z\langle x,y \rangle)/(N/Z\langle x,y \rangle)
\cong \langle x,y \rangle / N.
\]
The group 
\[
\langle x,y \rangle / Z\langle x,y \rangle \cong 
\langle \psi(x),\psi(y) \rangle
\]
has no non-trivial normal subgroups of infinite index by
Lemma~\ref{Lemma10}(2), hence $N/Z\langle x,y \rangle = 1$.
\end{itemize}
\end{proof}

\section{The example $\langle 1 + j + k, 1 + 2j \rangle$}
\label{SectionEx}
We study in this section the group $\langle 1 + j + k, 1 + 2j \rangle$.
Let $p = |1 + j + k|^2 = 3$ and $l = |1 + 2j|^2 = 5$. 
Then the group $\Gamma := \Gamma_{3,5}$
has the finite presentation
\begin{align}
\Gamma = \langle a_1, a_2, b_1, b_2, b_3 \mid 
&a_1 b_1 a_2 b_2, \, a_1 b_2 a_2 b_1^{-1}, \, a_1 b_3 a_2^{-1} b_1, \notag \\
&a_1 b_3^{-1} a_1 b_2^{-1}, \, a_1 b_1^{-1} a_2^{-1} b_3, \, a_2 b_3 a_2 b_2^{-1} \rangle, \notag
\end{align}
where we take
\begin{align}
a_1 &:= \psi_{3,5}( 1 + j + k), & a_1^{-1} &= \psi_{3,5}( 1 - j - k), \notag \\    
a_2 &:= \psi_{3,5}( 1 + j - k), & a_2^{-1} &= \psi_{3,5}( 1 - j + k), \notag \\
b_1 &:= \psi_{3,5}( 1 + 2i), & b_1^{-1} &= \psi_{3,5}( 1 - 2i), \notag \\    
b_2 &:= \psi_{3,5}( 1 + 2j), & b_2^{-1} &= \psi_{3,5}( 1 - 2j), \notag \\  
b_3 &:= \psi_{3,5}( 1 + 2k), & b_3^{-1} &= \psi_{3,5}( 1 - 2k).\notag
\end{align}
In the following, let $x := 1 + j + k$, $y := 1 + 2j$,
$a := a_1 = \psi_{3,5}(x)$, $b := b_2 = \psi_{3,5}(y)$ and
define the five words
\begin{align}
r_1(a,b) &:= ba^2bab^{-1}a^4b^{-1}a, \notag \\
r_2(a,b) &:= a^{-1}ba^{-1}ba^2ba^{-2}ba^{-1}b^2a^2bab, \notag \\
r_3(a,b) &:= baba^2b^2ab^{-1}ab^2a^2bab^2, \notag \\
r_4(a,b) &:= ba^2ba^{-1}b^{-3}a^{-2}b^{-1}ab^2, \notag \\
r_5(a,b) &:= ab^{-1}a^2b^{-1}ab^{-1}a^{-2}b^{-1}a^{-2}ba^{-2}ba^{-1}ba^2ba \notag
\end{align}
of lengths $12$, $18$, $18$, $14$ and $22$, respectively.
We will simply write $r_1, \ldots, r_5$ instead of 
$r_1(a,b), \ldots, r_5(a,b)$
or $r_1(x,y), \ldots, r_5(x,y)$ if the context is unambiguous,
like in the presentations of $\langle a,b \rangle$ and
$\langle x,y \rangle$ given below.
Using \textsf{GAP} (\cite{GAP}), we have done the following
computations.

\begin{Observation} \label{Obs12}
Let $\Gamma, a, a_2, b_1, b, b_3, x, y, r_1, r_2, r_3, r_4, r_5$ be as above. Then
\begin{itemize}
\item[(1)] $\Gamma^{ab} \cong \mathbb{Z}_2 \times \mathbb{Z}_4 \times \mathbb{Z}_4$.
\item[(2)] $(\Gamma')^{ab} \cong \mathbb{Z}_8 \times \mathbb{Z}_8 \times \mathbb{Z}_{16}$.
\item[(3)] $[\Gamma : \langle a,b \rangle] = 2$ such that 
$a_2, b_1 \notin \langle a,b \rangle$ and $b_3 = ab^{-1}a$.
\item[(4)] $\langle a,b \rangle$ has the finite presentation 
$\langle a, b \mid r_1, r_2, r_3, r_4, r_5 \rangle$.
There is no shorter non-trivial freely reduced relator than $r_1$.
\item[(5)] $\langle a,b \rangle^{ab} \cong 
\langle a, b \mid r_1, r_2, r_3, r_4, r_5, [a,b] \rangle
\cong \langle a, b \mid a^8, b^8, [a,b] \rangle
\cong \mathbb{Z}_8 \times \mathbb{Z}_8$.
\item[(6)] $\langle a, b \rangle' = \langle [a,b], [a,b^{-1}], [a^{-1},b], [a,b^2], [a^2,b] \rangle$.
\item[(7)] $(\langle a, b \rangle')^{ab} \cong \mathbb{Z}_{8} \times \mathbb{Z}_{8} \times \mathbb{Z}_{64}$.
\item[(8)] $r_1(x,y) = 3^4$, $r_2(x,y) = 5^4$, $r_3(x,y) = 3^{4} 5^{4}$, $r_4(x,y) = r_5(x,y) = 1$.
\end{itemize}
\end{Observation}
By Lemma~\ref{Lemma3}(3),(8) we have
\[
\mathbb{Z} \times \mathbb{Z} \cong \langle 3^4, 5^4 \rangle < Z\langle x,y\rangle < \langle -1,3,5 \rangle
\cong \mathbb{Z}_2 \times \mathbb{Z} \times \mathbb{Z}.
\]
It will turn out that 
\[
Z\langle x,y\rangle = \langle 3^4, 5^4 \rangle = \langle r_1(x,y),
r_2(x,y) \rangle.
\]
This enables us to compute an explicit presentation of the group 
$\langle x,y \rangle$ in the following proposition:

\begin{Proposition} \label{Prop13}
Let $\Gamma, a_1, a_2, b_1, b_2, b_3, x, y, a, b, r_1, r_2, r_3, r_4, r_5$ be as above. Then 
\begin{itemize}
\item[(1)] $Z\langle x,y\rangle = \langle 3^4, 5^4 \rangle$.
\item[(2)] $\langle x, y \rangle$ has a finite presentation
\[
\langle x, y \mid  r_4, \, r_5, \, r_1 r_2 r_3^{-1}, \, [x,r_1], \,
[x,r_2], \, [y,r_1], \, [y,r_2] \rangle
\]
\begin{align}
= \langle &x, y \mid \notag \\
&y x^2 y X Y^3 X^2 Y x y^2, \notag \\
&x Y x^2 Y x Y X^2 Y X^2 y X^2 y X y x^2 y x, \notag \\
&y x^2 y x Y x^3 y x^2 y X^2 y X y^2 x^2 y x Y
X Y X^2 Y^2 X y X Y^2 X^2 Y X Y,
\notag \\
&x y x^2 y x Y x^4 Y X y X^4 y X Y X^2 Y, \notag \\
&y X y x^2 y X^2 y X y^2 x^2 y x y X Y X
Y X^2 Y^2 x Y x^2 Y X^2 Y x Y x,
\notag \\
&y^2 x^2 y x Y x^4 Y x Y X y X^4 y X Y X^2 Y,
\notag \\
&y X y X y x^2 y X^2 y X y^2 x^2 y x Y X Y X^2
Y^2 x Y x^2 Y X^2 Y x Y x 
\rangle, \notag
\end{align}
where we write $X := x^{-1}$ and $Y := y^{-1}$. 
There is no shorter non-trivial freely reduced relator than $r_4(x,y)$ in $\langle x, y \rangle$.
\item[(3)] $\langle x, y \rangle$ is torsion-free.
\item[(4)] $\langle x, y \rangle^{ab} \cong \mathbb{Z} \times \mathbb{Z}$
and $(\langle x, y \rangle')^{ab} \cong \mathbb{Z}_8 \times
\mathbb{Z}_8 \times \mathbb{Z}_{64}$.
\item[(5)] $\langle x, y \rangle'$ has amalgam decompositions 
\[
F_{65} \ast_{F_{385}} F_{65} \text{\, and \,}
F_{129} \ast_{F_{513}} F_{129}.
\]
Moreover, the group $\langle x, y \rangle''$ has amalgam decompositions 
\[
F_{262145} \ast_{F_{1572865}} F_{262145} \text{\, and \,}
F_{524289} \ast_{F_{2097153}} F_{524289}.
\] 
\item[(6)] The two quaternions 
\[
y = 1+2j \text{\, and \,} xy^{-1}x = 3 \cdot 5^{-1}(1+2k)
\] 
generate a free subgroup $F_2$ of $\langle x, y \rangle$.
\item[(7a)] Let $r := [x, yx^{-1}y]$ and $q := x^{-1}rx$.
Then the two quaternions 
\[
r^2 q r^4 = 3^{-2} 5^{-12}(1700294841 + 519258632 i - 556215472 j + 1165319056 k)
\]
and 
\[
r^4 q r^2 = 3^{-2} 5^{-12}(1700294841 + 1191258632 i + 283784528 j + 661319056 k)
\]
generate a free subgroup $F_2$ of $\langle x, y \rangle'$.
\item[(7b)] The two quaternions
\[
[xy^{-1}x,y] = 5^{-2}(-7 - 8 i - 16 j + 16 k)
\]
and
\[
[y^{-1}, xy^{-1}x] = 5^{-2}(-7 - 8 i - 16 j - 16 k) 
\]
generate a free subgroup $F_2$ of $\langle x, y \rangle'$.
\item[(7c)] The two quaternions
\[
\frac{y^8}{5^4}
= 5^{-4}(-527 + 336 j)
\]
and
\[
\frac{5^4 (x y^{-1} x)^8}{3^8} 
= 5^{-4}(-527 + 336 k)
\]
generate a free subgroup $F_2$ of $\langle x, y \rangle'$.
\item[(8)] Let 
\[
w_1 := x y^{-1} x^{-1} y^{-1} x^{-2} y^{-2} = 3^{-3}5^{-2}(5 + 4i + 6j - 2k) 
\]
and
\[
w_2 := y^{-1} x^2 y^{-1} x y^{-1} x y^{-1} x^2 = -3^4 5^{-4}(11/3 + 4i + 6j - 2k).
\]
Then $\langle w_1, w_2 \rangle$ is a subgroup of $\langle x,y \rangle$
isomorphic to $\mathbb{Z} \times \mathbb{Z}$. 
\end{itemize}
\end{Proposition}

\begin{proof}
\begin{itemize}
\item[(1)]
Let $z := 1+j-k$, $s := 1+2i$, $t := 1+2k$,
and let 
\[
G := \langle x,z,s,y,t \rangle < Q_{3,5}.
\]
Our first goal is to obtain an explicit finite presentation of $G$.
This will be useful to compute $Z\langle x,y\rangle$.

Let $u_1, \ldots, u_6$ be the six defining relators of
$\Gamma$ from above, i.e.\
\begin{align}
u_1(a_1,a_2,b_1,b_2,b_3) &:= a_1 b_1 a_2 b_2, \notag \\
u_2(a_1,a_2,b_1,b_2,b_3) &:= a_1 b_2 a_2 b_1^{-1}, \notag \\
u_3(a_1,a_2,b_1,b_2,b_3) &:= a_1 b_3 a_2^{-1} b_1, \notag \\
u_4(a_1,a_2,b_1,b_2,b_3) &:= a_1 b_3^{-1} a_1 b_2^{-1}, \notag \\
u_5(a_1,a_2,b_1,b_2,b_3) &:= a_1 b_1^{-1} a_2^{-1} b_3, \notag \\
u_6(a_1,a_2,b_1,b_2,b_3) &:= a_2 b_3 a_2 b_2^{-1}. \notag
\end{align}
As in Lemma~\ref{Lemma3}(3), it is easy to see that $ZG < \langle -1, 3, 5 \rangle$.
On the other hand, we check that
\begin{align}
-1 &= u_5(x,z,s,y,t) = x s^{-1} z^{-1} t \in G, \notag \\
 3 &= u_6(x,z,s,y,t) = z t z y^{-1} \in G, \notag \\
 5 &= -u_3(x,z,s,y,t) = -x t z^{-1} s \in G. \notag
\end{align}
This shows that $ZG = \langle -1, 3, 5 \rangle$. 
We note that this implies $G = Q_{3,5}$ and that $G$ is not torsion-free.

Let $\psi$ be the restriction of $\psi_{3,5}$ to $G$. 
As in Lemma~\ref{Lemma3}(2),
we have 
\[
\mathrm{ker}(\psi : G \twoheadrightarrow \Gamma) = ZG.
\]
Knowing finite presentations of $\Gamma$ and $ZG$, it is now possible to compute
a finite presentation of the extension $G$ 
(for example see \cite[Chapter~10.2]{Johnson} 
for the detailed explicit general construction
of such a presentation).
Evaluating the three other defining relators of $\Gamma$
\begin{align}
u_1(x,z,s,y,t) &= x s z y = -15, \notag \\
u_2(x,z,s,y,t) &= x y z s^{-1} = -3, \notag \\
u_4(x,z,s,y,t) &= x t^{-1} x y^{-1} = 3/5, \notag
\end{align}
we get the finite presentation
\begin{align}
G = \langle x,z,s,y,t \mid \, &[x,u_3], \, [x,u_5], \, [x,u_6], \notag
                           \\ &[z,u_3], \, [z,u_5], \, [z,u_6], \notag \\
                           &[s,u_3], \, [s,u_5], \, [s,u_6], \notag\\ &[y,u_3], \, [y,u_5], \, [y,u_6], \notag \\
                           &[t,u_3], \, [t,u_5], \, [t,u_6], \notag \\ &u_5^2, \,
                           u_1u_3^{-1}u_6^{-1}, \, u_2u_5u_6^{-1}, \, u_3u_4u_5u_6^{-1} \rangle. \notag
\end{align}

We have seen that $\langle 3^4, 5^4 \rangle < Z\langle x,y\rangle < \langle -1,3,5 \rangle$.
The group $\langle 3^4, 5^4 \rangle$ has index $32 = 2 \cdot 4 \cdot 4$ in $\langle -1,3,5 \rangle$,
and we write $\langle -1,3,5 \rangle$ as a finite disjoint union of cosets
\[  
\langle -1,3,5 \rangle = \langle 3^4, 5^4 \rangle \sqcup \lambda_2 \langle 3^4, 5^4 \rangle \sqcup \ldots \sqcup
\lambda_{32} \langle 3^4, 5^4 \rangle,
\]
such that $\lambda_2, \ldots, \lambda_{32} \in \langle -1,3,5 \rangle
\setminus \langle 3^4, 5^4 \rangle$.
To prove our statement $Z \langle x, y \rangle = \langle 3^4, 5^4 \rangle$,
it is enough to check that the index condition
\[
[G : \langle x, y \rangle ] > [G : \langle x, y, \lambda_m \rangle]
\]
is satisfied for all $m \in \{2, \ldots, 32\}$.
To see why this is enough, assume that this index condition holds
and that $Z \langle x, y \rangle \ne \langle 3^4, 5^4 \rangle$.
Then there is an element $\gamma \in Z \langle x, y \rangle \setminus \langle 3^4, 5^4 \rangle$,
hence $\gamma \in \lambda_m \langle 3^4, 5^4 \rangle$ for some $m \in \{2, \ldots, 32\}$,
and therefore $\lambda_m \in Z \langle x, y \rangle < \langle x, y \rangle$,
which is impossible, since the index condition implies that 
$\lambda_m \notin \langle x, y \rangle$ for all $m \in \{2, \ldots, 32\}$.

Using \textsf{GAP} (\cite{GAP}) and the finite presentation of $G$ from above, 
we have checked that the condition
$[G : \langle x, y \rangle ] > [G : \langle x, y, \lambda_m \rangle]$ indeed holds.
In particular, we have $[G : \langle x, y \rangle] = 64$
and $[G : \langle x, y, \lambda \rangle] = 32 < 64$,
if
\[
\lambda \in \{ -1, 3^2, -3^2, 5^2, -5^2, 3^25^2, -3^25^2\}.
\]
\item[(2)]
We obtain the stated finite presentation of $\langle x, y \rangle$
as an extension of $Z\langle x, y \rangle$ by $\langle a,b \rangle$.
Observe that the four commutator relators in the presentation express
the fact that $r_1 = 3^4$ and $r_2 = 5^4$ are in the center of $\langle x, y \rangle$.

We have checked with \textsf{GAP} (\cite{GAP}) that there is no 
non-trivial freely reduced relator of length less than $14$. 
Moreover, any freely reduced relator of length $14$
is conjugate to $r_4$ or to $r_4^{-1}$.
\item[(3)]
The group $\langle x, y \rangle$ is torsion-free by Lemma~\ref{Lemma8}, 
using that $-1 \notin Z\langle x, y \rangle$.
\item[(4)]
The abelianization of $\langle x, y \rangle$ is easy to compute, once the finite presentation
of $\langle x, y \rangle$ given in part~(2) of this proposition is known.
Alternatively, we apply Proposition~\ref{Prop9}(1).

Since $-1 \notin Z\langle x, y \rangle$, we have
\[ 
\langle x,y \rangle' \cap Z\langle x, y \rangle = 1,
\]
hence an isomorphism $\langle a,b \rangle' \cong \langle x,y \rangle'$
by the left column of our diagram of Lemma~\ref{Lemma3}(5).
Now, we use Observation~\ref{Obs12}(7).
\item[(5)]
Since $\langle x,y \rangle' \cong \langle a,b \rangle'$, 
it is enough to show that $\langle a,b \rangle'$ has the claimed amalgam decompositions.
Let $\Gamma_0$ be the kernel of the surjective homomorphism
$\Gamma \to \mathbb{Z}_2 \times \mathbb{Z}_2$, given by
\[
a_1, a_2 \mapsto (1 + 2\mathbb{Z}, 0 + 2\mathbb{Z}) \text{\, and \,} 
b_1, b_2, b_3 \mapsto (0 + 2\mathbb{Z}, 1 + 2\mathbb{Z}).
\]
Then $\Gamma_0$ has by \cite[Proposition~1.4]{Rattaggi} amalgam decompositions
$F_{3} \ast_{F_{13}} F_{3}$ and
$F_{5} \ast_{F_{17}} F_{5}$.
Commutators in $\Gamma$ are obviously in $\Gamma_0$,
in particular $\langle a,b \rangle' < \Gamma_0$.
We have checked with \textsf{GAP} (\cite{GAP}) that 
$\langle a,b \rangle' = \langle \! \langle [a,b] \rangle \!
\rangle_{\langle a,b \rangle}$
is in fact $\langle \! \langle [a,b] \rangle \! \rangle_{\Gamma}$,
the normal closure of $[a,b]$ in $\Gamma$.
It follows that $\langle a,b \rangle'$ is a normal subgroup of
$\Gamma$, hence a normal subgroup of $\Gamma_0$.
The index of $\langle a,b \rangle'$ in $\Gamma_0$ is
\[
[\Gamma_0 : \langle a,b \rangle'] = \frac{[\Gamma : \langle a,b \rangle']}{[\Gamma : \Gamma_0]}
= \frac{[\Gamma : \langle a,b \rangle] \cdot |\langle a,b \rangle^{ab}|}{[\Gamma : \Gamma_0]}
= \frac{2 \cdot 64}{4} = 32.
\]
Now, we apply \cite[Corollary~2, Corollary~1]{Djokovic} to get the stated amalgam decompositions
of $\langle a,b \rangle' \cong \langle x,y \rangle'$, observing that
\[
32 = [F_{3} : F_{65}] = [F_{13} : F_{385}] = [F_{5} : F_{129}] = [F_{17} : F_{513}].
\]

The second claim follows similarly, using that $\langle x,y \rangle''$
has index $4096 = 8 \cdot 8 \cdot 64$ in $\langle x,y \rangle'$ 
by part~(4) of this proposition.
\item[(6)]
A relation in $y$ and $xy^{-1}x$ induces a relation in 
$\psi(y) = b_2$ and $\psi(xy^{-1}x) = a_1 b_2^{-1} a_1 = b_3$, but 
$F_2 \cong \langle b_1, b_2 \rangle < \langle b_1, b_2, b_3
\rangle \cong F_3 < \Gamma$.
\item[(7a)]
We compute
\[
r = [x,yx^{-1}y] = - \frac{7}{25} + \frac{8}{75}i + \frac{32}{75}j + \frac{64}{75}k \in \langle x,y \rangle'
\]
and
\[
q = x^{-1}rx = 
[y,x^{-1}yx^{-1}] = - \frac{7}{25} - \frac{8}{25}i + \frac{16}{25}j + \frac{16}{25}k \in \langle x,y \rangle',
\]
clearly both of norm $1$.
Since 
\[
r = - \frac{7}{25} + \frac{24}{25}\left( \frac{1}{9}i + \frac{4}{9}j + \frac{8}{9}k\right)
\]
and
\[
q = - \frac{7}{25} - \frac{24}{25}\left( \frac{1}{3}i - \frac{2}{3}j - \frac{2}{3}k\right),
\]
they are both ``rational'' in the sense of \cite[Definition~4.1]{Daless}, and by
\cite[Corollary~4.2]{Daless} we have
\[
F_2 \cong \langle r^2 q r^4, r^4 q r^2 \rangle < \langle x,y \rangle' < \langle x,y \rangle.
\]
\item[(7b)]
We have $\psi([xy^{-1}x,y]) = b_3 b_2 b_3^{-1} b_2^{-1}$
and $\psi([y^{-1}, xy^{-1}x]) = b_2^{-1} b_3 b_2 b_3^{-1}$.
The group
$\langle b_3 b_2 b_3^{-1} b_2^{-1}, b_2^{-1} b_3 b_2 b_3^{-1} \rangle$
is free as a subgroup of $\langle b_2, b_3 \rangle \cong F_2$, 
but not isomorphic to $\mathbb{Z}$,
hence 
\[
\langle b_3 b_2 b_3^{-1} b_2^{-1}, b_2^{-1} b_3 b_2 b_3^{-1} \rangle \cong F_2.
\]
By Lemma~\ref{Lemma3}(9)
\[
\langle [xy^{-1}x,y], [y^{-1}, xy^{-1}x] \rangle \cong F_2.
\]
This is clearly a subgroup of $\langle x,y \rangle'$.
\item[(7c)]
Let 
\[
q_1 := \frac{y^8}{5^4}
\text{\, and \,}
q_2 := \frac{5^4 (x y^{-1} x)^8}{3^8}.
\]
First recall that $3^4, 5^4 \in \langle x,y \rangle$.
It follows that $q_1, q_2 \in \langle x, y \rangle$.
We have
\[
\psi(q_1) = \psi(y^8) = b_2^8 
\]
and
\[
\psi(q_2) = \psi((x y^{-1} x)^8) = b_3^8. 
\]
Since $\langle b_2^8, b_3^8 \rangle \cong F_2$,
Lemma~\ref{Lemma3}(9) implies
$\langle q_1, q_2 \rangle \cong F_2$.
It is easy to see that $q_1$ and $q_2$ both have norm $1$.
In particular $\langle q_1, q_2 \rangle$ is a free subgroup
of $\langle x, y \rangle \cap \mathbb{H}(\mathbb{Q})_1$,
i.e.\ a free subgroup of $\langle x, y \rangle'$
using Proposition~\ref{Prop9}(2).
\item[(8)] This example is an illustration of the proof of
Theorem~\ref{Thm11}(6). In \cite[Section~4.1]{Rattaggi} we have shown
that 
\[
\Gamma > \langle a_1 a_2 a_1 a_2^{-1}, b_2^{-1} b_1^{-1} b_3^{-1} b_1 \rangle
\cong \mathbb{Z} \times \mathbb{Z}.
\]
Using \textsf{GAP} (\cite{GAP}) we check that
\[
[\Gamma : \langle a, b, a_1 a_2 a_1 a_2^{-1} \rangle] = 
[\Gamma : \langle a, b, b_2^{-1} b_1^{-1} b_3^{-1} b_1 \rangle] 
= [\Gamma : \langle a, b \rangle] = 2,
\] 
hence $a_1 a_2 a_1 a_2^{-1} \in \langle a, b \rangle$ and
$b_2^{-1} b_1^{-1} b_3^{-1} b_1 \in \langle a, b \rangle$.
Indeed, it is easy to check that
\[
\psi(w_1) = a b^{-1} a^{-1} b^{-1} a^{-2} b^{-2} = a_1 a_2 a_1 a_2^{-1}
\]
and
\[
\psi(w_2) = b^{-1} a^2 b^{-1} a b^{-1} a b^{-1} a^2 = b_2^{-1} b_1^{-1} b_3^{-1} b_1.
\]
\end{itemize}
\end{proof}

\begin{Remark} \label{Rem14}
Using the presentation of $G = Q_{3,5}$ computed in the proof of
Proposition~\ref{Prop13}(1), 
it is easy to check that $G^{ab}$ and 
$(G/\langle \! \langle -1 \rangle \! \rangle_G)^{ab}$
are not isomorphic. 
Indeed, 
$G^{ab} \cong \mathbb{Z} \times \mathbb{Z} \times \mathbb{Z}_2 \times \mathbb{Z}_2 \times \mathbb{Z}_4$,
whereas 
$(G/\langle \! \langle -1 \rangle \! \rangle_G)^{ab} \cong 
\mathbb{Z} \times \mathbb{Z} \times \mathbb{Z}_2 \times \mathbb{Z}_4$.
This shows that $-1 \notin G' = {Q_{3,5}}'$.
In the same way, we have also checked that 
$-1 \notin {Q_{p,l}}'$, if $(p,l) = (3,7)$, $(3,11)$, $(5,7)$, or $(5,11)$,
in particular Conjecture~\ref{Conj5} and Conjecture~\ref{Conj6} 
are true for these pairs $(p,l)$.
\end{Remark}

\newpage

See the big diagram below for the relations of some subgroups of
$\Gamma$ and $Q_{3,5}$ used in some previous proofs.
In this diagram, all simple arrows are injective homomorphisms.
A label of the form $[n_1, \ldots, n_k]$ stands for 
the quotient group 
$\mathbb{Z}/n_1 \mathbb{Z} \times \ldots \times \mathbb{Z}/n_k \mathbb{Z}$.
The other labels of the injective homomorphisms are indices.
The labels of the two surjective homomorphisms on top right indicate
the corresponding kernels.
The powers of $2$ on the left are the indices of the subgroups of
$\Gamma$ on the same line. 
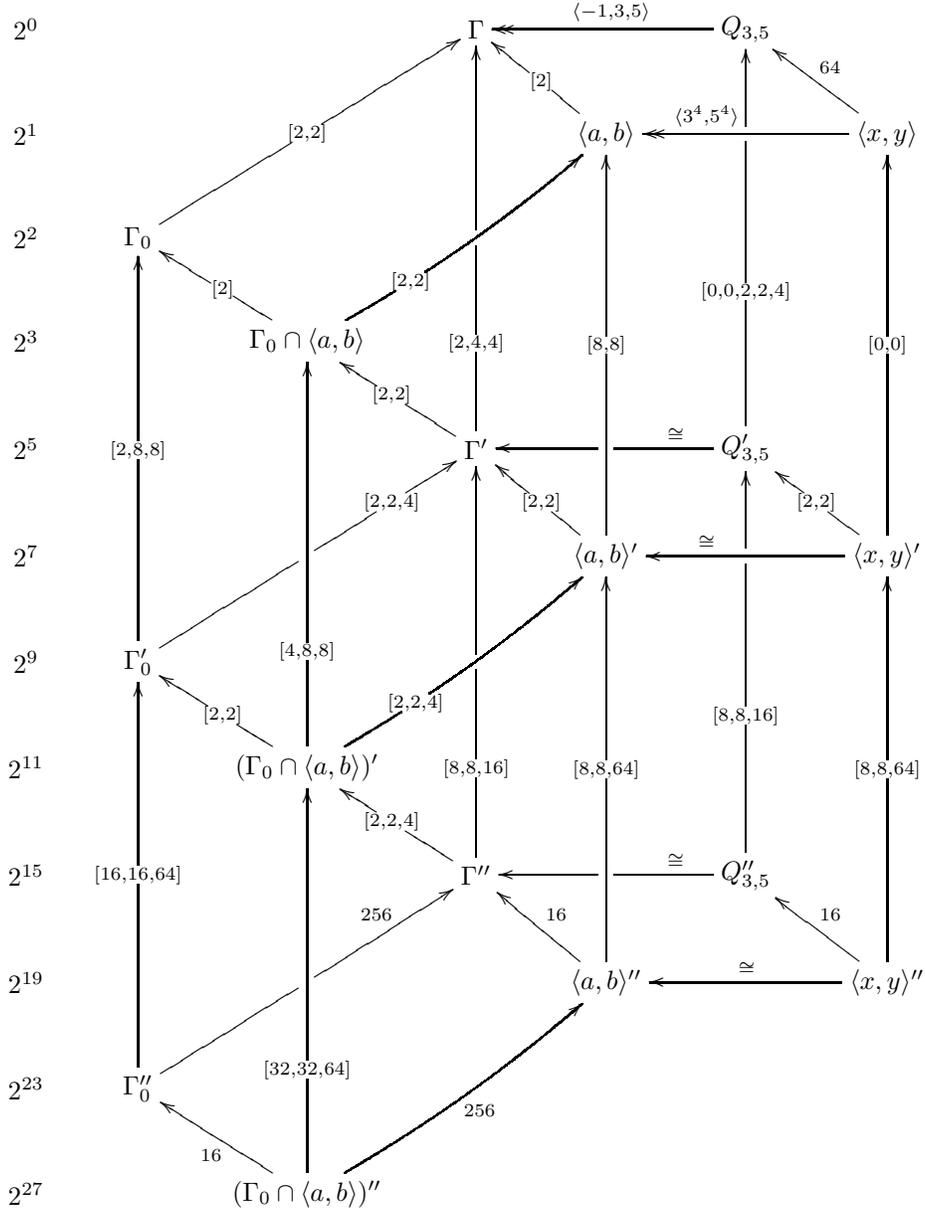
\begin{figure}[ht]
\[
\xymatrix{
2^0&&& \Gamma&& Q_{3,5} \ar@{->>}[ll]_{\langle -1,3,5 \rangle}&\\
2^1 &&&& \langle a,b \rangle \ar[lu]|{[2]} && \langle x,y
\rangle \ar@{->>}[ll]_{\langle 3^4,5^4 \rangle \phantom{aaaaaaa}}\ar[ul]_{64} &\\
2^2 & \Gamma_0 \ar[rruu]|{[2,2]}&&&&&\\
2^3 && \Gamma_0 \cap \langle a,b \rangle \ar[lu]|{[2]} \ar@/_0.35pc/[rruu]|<<<<<<<{[2,2]} &&&&\\
2^5 &&& \Gamma' \ar'[uu]|{[2,4,4]}[uuuu] \ar[lu]|{[2,2]} && Q_{3,5}'
\ar'[l]_{\cong}[ll] \ar'[uuu]|{[0,0,2,2,4]}[uuuu]&\\
2^7 &&&& \langle a,b \rangle' \ar[uuuu]|{[8,8]}
\ar[lu]|{[2,2]} && \langle x,y \rangle'\ar[ll]_{\cong \phantom{aaaaaaa}} \ar[lu]|{[2,2]} \ar[uuuu]|{[0,0]}&\\
2^9    & \Gamma_0' \ar[uuuu]|{[2,8,8]} \ar'[ur][rruu]|{[2,2,4]}  &&&&&\\
2^{11} && (\Gamma_0 \cap \langle a,b \rangle)'
\ar[uuuu]|<<<<<<<<<<<<<{[4,8,8]} \ar[lu]|{[2,2]} \ar@/_0.34pc/[rruu]|<<<<<{[2,2,4]}&&&&\\
2^{15} &&& \Gamma''  \ar'[uu]|{[8,8,16]}[uuuu] \ar[lu]|{[2,2,4]} && Q_{3,5}''
\ar'[l]_{\cong}[ll] \ar'[uuu]|{[8,8,16]}[uuuu]&\\
2^{19} &&&& \langle a,b \rangle'' \ar[uuuu]|{[8,8,64]} \ar[lu]_{16}&& \langle x,y \rangle''
\ar[ll]_{\cong} \ar[lu]_{16} \ar[uuuu]|{[8,8,64]}&\\
2^{23} & \Gamma_0'' \ar[uuuu]|{[16,16,64]} \ar'[ru][rruu]^{256} &&&&&\\
2^{27} && (\Gamma_0 \cap \langle a,b \rangle)''
\ar[lu]^{16} \ar@/_0.34pc/[rruu]_{256} \ar[uuuu]|<<<<<<<<<<<<<<{[32,32,64]}&&&&
}
\]
\end{figure}

\newpage

\section{The commutative case} \label{SectionComm}
In this section, we study the following simple question:
\begin{Question} \label{Question15}
Given two distinct odd prime numbers $p$ and $l$, 
are there commuting quaternions $x \in X_p$, $y \in X_l$?
\end{Question}

In the following, we will give general (negative and positive) 
answers to Question~\ref{Question15}, except in the three cases
$p \equiv 1, l \equiv 7 \pmod{8}$,
$p \equiv 7, l \equiv 1 \pmod{8}$
and $p, l \equiv 7 \pmod{8}$,
where the situation seems to be quite complicated.

Let $T_{p,l}$ be the set $Q_{p,l} \cap \mathbb{H}(\mathbb{Z})$, i.e.\
\begin{align*}
T_{p,l} := \{x = x_0  + &x_1 i + x_2 j + x_3 k \in \mathbb{H}(\mathbb{Z}) \,; 
\quad |x|^2 = p^r l^s, \; r,s \in \mathbb{N}_0; \notag \\
&x_0 \text{ odd}, x_1, x_2, x_3 \text{ even}, \text{ if } |x|^2 \equiv 1 \!\!\!\! \pmod 4\,; \notag \\
&x_1 \text{ even}, x_0, x_2, x_3 \text{ odd}, \text{ if } |x|^2 \equiv 3 \!\!\!\! \pmod 4
\} \, .
\end{align*}
Then clearly $X_p \subset T_{p,l}$ and $X_l \subset T_{p,l}$.
Note that $x,y \in T_{p,l}$ 
commute, if and only if $\psi_{p,l}(x), \psi_{p,l}(y)$ commute in $\Gamma_{p,l}$,
but we will directly work with quaternions here and not use the group $\Gamma_{p,l}$ anymore.

Let $x = x_0 + x_1 i + x_2 j + x_3 k \in T_{p,l}$ such that $(x_1,x_2,x_3) \ne (0,0,0)$.
We can write it as 
\[
x = x_0 + z_x(c_1 i + c_2 j + c_3 k)
\]
such that $c_1, c_2, c_3 \in \mathbb{Z}$ are
relatively prime and $z_x \in \mathbb{Z} \setminus \{ 0 \}$.
Up to multiplication by $-1$, the integers
$c_1, c_2, c_3, z_x$ are uniquely determined by $x$.
Therefore we can define the number 
\[
n(x) := c_1^2 + c_2^2 + c_3^2 \in \mathbb{N}.
\]
Moreover, if 
$x = x_0 \in T_{p,l} \cap \mathbb{Q}^{\ast} = T_{p,l} \cap \mathbb{Z}$, 
we define $n(x) := 0$.

\begin{Remark} \label{Rem16}
It is easy to show that $n(x) \equiv 1,2,3,5,6 \pmod{8}$, if $n(x) \ne 0$. 
\end{Remark}
The function $n : T_{p,l} \to \mathbb{N}_0$
is a useful invariant for commutativity:

\begin{Lemma} \label{Lemma17}
Let $p,l$ be two distinct odd prime numbers and
let $x,y \in T_{p,l}$ be commuting quaternions such that $x, y \notin \mathbb{Q}^{\ast}$. 
Then $n(x) = n(y) \ne 0$.
\end{Lemma} 

\begin{proof}
This follows directly from the basic fact 
(see \cite[Section~3]{Mozes1} or \cite[Lemma~12]{Rattaggi2})
that two quaternions 
$x = x_0 + x_1 i + x_2 j + x_3 k$ and
$y = y_0 + y_1 i + y_2 j + y_3 k$
in $\mathbb{H}(\mathbb{Q}) \setminus \mathbb{Q}$ commute, if and only if
$\mathbb{Q}(x_1,x_2,x_3) = \mathbb{Q}(y_1,y_2,y_3)$.
\end{proof}

Lemma~\ref{Lemma17} will be used in Proposition~\ref{Prop21} to show that
certain pairs $x \in X_p$, $y \in X_l$ cannot commute. To apply this lemma we have to
get some knowledge on~$n$. This is
done in Lemma~\ref{Lemma19} for the cases where $|x|^2 \not\equiv 1 \pmod{8}$.
First we state an auxiliary very basic lemma:

\begin{Lemma} \label{Lemma18}
If $x_0 \in \mathbb{Z}$ is odd, then $x_0^2 \equiv 1 \pmod{8}$.
\end{Lemma}

\begin{proof}
Let $x_0 = 1+2t$ for some $t \in \mathbb{Z}$.
Then $x_0^2 = 1+4t(1+t) \equiv 1 \pmod{8}$, since $t(1+t)$ is always even.
\end{proof}

\begin{Lemma} \label{Lemma19}
Let $p$ be an odd prime number 
and let $x \in X_p$.
\begin{itemize}
\item[(1)] If $p \equiv 5 \pmod{8}$, then $n(x)$ is odd.
\item[(2)] If $p \equiv 3 \pmod{8}$, then $n(x) \equiv 2 \pmod{8}$.
\item[(3)] If $p \equiv 7 \pmod{8}$, then $n(x) \equiv 6 \pmod{8}$.
\end{itemize}
\end{Lemma}

\begin{proof}
First note that $X_p \cap \mathbb{Q} = \emptyset$.
Write $x = x_0 + z_x(c_1 i + c_2 j + c_3 k)$ as in the definition of $n$.
\begin{itemize}
\item[(1)] Since $|x|^2 \equiv 1 \pmod{4}$, $x_0$ is odd and $z_x$ is even.
We have 
\[
p = |x|^2 = x_0^2 + z_x^2n(x) \equiv 5 \pmod{8}.
\]
If $n(x)$ would be even, then $z_x^2n(x) \equiv 0 \pmod{8}$, 
hence $x_0^2 \equiv 5 \pmod{8}$, contradicting Lemma~\ref{Lemma18}.
\item[(2)] Here, $|x|^2 \equiv 3 \pmod{4}$, hence $x_0$, $z_x$ are odd and $n(x)$ is even 
(as a sum of two odd squares and one even square).
We have 
\[
p = |x|^2 = x_0^2 + z_x^2n(x) \equiv 3 \pmod{8},
\]
hence $z_x^2n(x) \equiv 2 \pmod{8}$ using Lemma~\ref{Lemma18}.
Since $z_x^2 \equiv 1 \pmod{8}$ by Lemma~\ref{Lemma18}, it follows that $n(x) \equiv 2 \pmod{8}$.
\item[(3)]  The proof is completely analogous to the proof of part~(2).
\end{itemize}
\end{proof}

\begin{Remark} \label{Rem20}
We will later see from Table~\ref{Table2b} that all of the
possibilities $n(x) \equiv 1,3,5 \pmod{8}$ 
(in view of Remark~\ref{Rem16} and Lemma~\ref{Lemma19}(1))
can be realized if $p \equiv 5 \pmod{8}$.
Moreover, in the case $p \equiv 1 \pmod{8}$ not treated in
Lemma~\ref{Lemma19}, all possibilities $n(x) \equiv 1,2,3,5,6 \pmod{8}$ can be realized.
\end{Remark}

\begin{Proposition} \label{Prop21}
Let $p,l$ be two distinct odd prime numbers. 
Suppose that
$p, l \not\equiv 1 \pmod{8}$ and $p \not\equiv l \pmod{8}$. Then 
there are no commuting quaternions $x \in X_p$, $y \in X_l$.
\end{Proposition}

\begin{proof}
This follows directly from Lemma~\ref{Lemma17}, using Lemma~\ref{Lemma19}.
\end{proof}

To obtain positive answers to Question~\ref{Question15}, we will use some known
results on prime numbers of the form $r^2 + m s^2$, first for $m=1$
and $m=2$ in Proposition~\ref{Prop23}, later for $m=6$ and $m=22$ in Proposition~\ref{Prop24}.

\begin{Lemma} \label{Lemma22} (Fermat, see \cite[(1.1)]{Cox})
Let $p$ be an odd prime number.
There are $x_0, x_1 \in \mathbb{Z}$ such that $x_0^2 + x_1^2 = p$, if and only if $p \equiv 1 \pmod{4}$.
There are $x_0, z \in \mathbb{Z}$ such that $x_0^2 + 2z^2 = p$, if and only if $p \equiv 1,3 \pmod{8}$.
\end{Lemma}

This lemma can be applied as follows:

\begin{Proposition} \label{Prop23}
Let $p,l$ be two distinct odd prime numbers. 
Suppose that either $p,l \equiv 1 \pmod{4}$
or that $p,l \equiv 1, 3 \pmod{8}$. Then
there are commuting quaternions $x \in X_p$, $y \in X_l$.
\end{Proposition}

\begin{proof}
If $p,l \equiv 1 \pmod{4}$ then by Lemma~\ref{Lemma22}, there are $x_0$, $y_0$ odd, $x_1$, $y_1$ even, 
such that 
\[
x_0^2 + x_1^2 = p \text{ \, and \, } 
y_0^2 + y_1^2 = l.
\] 
Now we take the commuting quaternions $x = x_0 + x_1 i \in X_p$ and $y = y_0 + y_1 i \in X_l$.

If $p \equiv 1 \pmod{8}$, then by Lemma~\ref{Lemma22} there are $x_0, z \in \mathbb{Z}$ such that $x_0^2 + 2z^2 = p$.
It follows that $x_0$ is odd, hence $2 z^2 \equiv 0 \pmod{8}$ by Lemma~\ref{Lemma18} 
and $z$ is even (but non-zero).
We choose 
\[
x := x_0 + z(j+k) \in X_p,
\] 
in particular $|x|^2 = x_0^2 + 2z^2 = p$ and $n(x) = 2$.

If $p \equiv 3 \pmod{8}$, then again $x_0^2 + 2z^2 = p$ by Lemma~\ref{Lemma22}, but here
$x_0$ and $z$ are both odd, and we take as above 
\[
x := x_0 + z(j+k) \in X_p,
\]
such that $|x|^2 = p$ and $n(x) = 2$.

In the same way we construct 
\[
y := y_0 + z_y (j+k) \in X_l
\]
such that $|y|^2 = l \equiv 1, 3 \pmod{8}$ and $n(y) = 2$. 
Clearly $xy = yx$ by construction.
\end{proof}

We illustrate the results of Proposition~\ref{Prop21} and Proposition~\ref{Prop23} in Table~\ref{Tablepl}
for $p,l$ taken modulo $8$
(``$+$'' means that there are always commuting quaternions $x \in X_p$, $y \in X_l$, 
``$-$'' means that there are never such commuting quaternions,
``$\pm$'' means that both cases happen).
\begin{table}[ht]
\[
\begin{tabular}{c|cccc}
$\pmod{8}$ & $l \equiv 1$ & $3$ & $5$ & $7$ \\ \hline
$p \equiv 1$              & $+$ & $+$ & $+$ & $\pm$ \\
$3$                       & $+$ & $+$ & $-$ & $-$  \\
$5$                       & $+$ & $-$ & $+$ & $-$  \\
$7$                       & $\pm$ & $-$ & $-$ & $\pm$ 
\end{tabular}
\]   
\caption{Existence and non-existence of commuting quaternions} \label{Tablepl}
\end{table}

In the three cases
$p \equiv 1, l \equiv 7 \pmod{8}$,
$p \equiv 7, l \equiv 1 \pmod{8}$
and  
$p, l \equiv 7 \pmod{8}$
excluded from Proposition~\ref{Prop21} and Proposition~\ref{Prop23},
it seems to be more difficult to decide in general by congruence conditions 
whether or not there are commuting quaternions $x \in X_p$, $y \in X_l$.
We illustrate this with results for some subcases of $p, l \equiv 7 \pmod{8}$
and $p \equiv 1, l \equiv 7 \pmod{8}$:

\begin{Proposition} \label{Prop24}
Let $p,l$ be two distinct odd prime numbers. 
\begin{itemize}
\item[(1)]
Suppose that
$p,l \equiv 7 \pmod{24}$ or 
$p,l \equiv 15, 23, 31, 47, 71 \pmod{88}$. Then
there are commuting quaternions $x \in X_p$, $y \in X_l$.
\item[(2)] Suppose that
$p \equiv 1, l \equiv 7 \pmod{24}$ or 
$p \equiv 1, 9, 25, 49, 81 \pmod{88}$, $l \equiv 15, 23, 31, 47, 71 \pmod{88}$. Then
there are commuting quaternions $x \in X_p$, $y \in X_l$.
\end{itemize}
\end{Proposition}

\begin{proof}
\begin{itemize}
\item[(1)]
A prime number $p$ is of the form $x_0^2 + 6z^2$, if and only if $p \equiv 1, 7 \pmod{24}$
(see \cite{Cox} for a proof, but beware of the misprint in \cite[(2.28)]{Cox}
stating the condition $\! \!  \pmod{12}$ instead of $\! \! \pmod{24}$).
Let $p,l \equiv 7 \pmod{24}$.
There are $x_0, z \in \mathbb{Z}$
such that 
\[
p = x_0^2 + 6z^2 =  x_0^2 + (2z)^2 + z^2 + z^2.
\]
It follows that $x_0$ and $z$ are odd.
Take 
\[
x := x_0 + z(2i + j + k) \in X_p,
\] 
hence $|x|^2 = x_0^2 + 6z^2 = p$ and $n(x) = 6$.
In the same way, we choose 
\[
y := y_0 + z_y(2i + j + k) \in X_l
\]
such that $|y|^2 = y_0^2 + 6z_y^2 = l$, $n(y) = 6$
and $xy = yx$.

If $p,l \equiv 15, 23, 31, 47, 71 \pmod{88}$, then we can give a similar proof,
using the fact (see \cite[(2.28)]{Cox}) that a prime number $p$ is of the form $x_0^2 + 22z^2$, if and only if 
$p \equiv 1, 9, 15, 23, 25, 31, 47, 49, 71, 81 \pmod{88}$.
Observe that $22z^2 \equiv 6 \pmod{8}$ if $z$ is odd, and
$22z^2 \equiv 0 \pmod{8}$ if $z$ is even.
Therefore $z$ is odd here.
We take 
\[
x := x_0 + z(2i + 3j + 3k) \in X_p,
\] 
hence $|x|^2 = x_0^2 + 22z^2 = p$ and $n(x) = 22$.
Now, we are done by the analogous construction for $y$.
\item[(2)]
Let $p \equiv 1 \pmod{24}$.
The proof is similar as in~(1), but here $p = x_0^2 + 6z^2$ for some
$x_0$ odd, $z =: 2\tilde{z}$ even (and non-zero),
hence 
\[
p = x_0^2 + 24 \tilde{z}^2 = x_0^2 + (4\tilde{z})^2 + (2\tilde{z})^2 +
(2\tilde{z})^2
\]
and we choose 
\[
x := x_0 + \tilde{z}(4i + 2j + 2k)
\]
which commutes with $2i+j+k$.

Let $p \equiv 1, 9, 25, 49, 81 \pmod{88}$. Similarly as above we can choose
\[
x := x_0 + \tilde{z}(4i+6j+6k)
\] 
of norm $x_0^2 + 22(2 \tilde{z})^2$
commuting with $2i +3j + 3k$.
\end{itemize}
\end{proof}

The two quaternions $x = 1+2i$ and $y = 1+4k$ do not commute, but $n(x) = n(y) = 1$, so the converse of
Lemma~\ref{Lemma17} is certainly not true. However, the function $n$ determines the \emph{existence} 
of commuting quaternions of given norms.

\begin{Proposition} \label{Prop25}
Let $p,l$ be two distinct odd prime numbers and
let $x,y \in T_{p,l}$ such that $n(x) = n(y)$.
Then there are commuting quaternions $\hat{x}, \hat{y} \in T_{p,l}$
such that $|\hat{x}|^2 = |x|^2$ and $|\hat{y}|^2 = |y|^2$.
\end{Proposition}

\begin{proof}
If $n(x)=n(y)=0$, then $x=x_0$, $y=y_0$
and we can take $\hat{x} := x$, $\hat{y} := y$.

Now suppose that $n(x)=n(y) \ne 0$.
Write 
\[
x = x_0 + x_1 i + x_2 j + x_3 k = x_0 + z_x(c_1 i + c_2 j + c_3 k)
\]
such that $c_1, c_2, c_3 \in \mathbb{Z}$ are
relatively prime and $z_x \in \mathbb{Z} \setminus \{ 0 \}$,
and write
\[
y = y_0 + y_1 i + y_2 j + y_3 k = y_0 + z_y(d_1 i + d_2 j + d_3 k)
\]
such that $d_1, d_2, d_3 \in \mathbb{Z}$ are
relatively prime and $z_y \in \mathbb{Z} \setminus \{ 0 \}$.
Then by assumption
\[
c_1^2 + c_2^2 + c_3^2 = n(x) = n(y) = d_1^2 + d_2^2 + d_3^2.
\]

If $|x|^2 \equiv 1 \pmod{4}$, then $z_x$ is even (and non-zero). Let
\[
\hat{x} := x_0 + z_x(d_1 i + d_2 j + d_3 k) \in T_{p,l}.
\]
Then 
\[
|\hat{x}|^2 = x_0^2 + z_x^2n(y) = x_0^2 + z_x^2n(x) = |x|^2
\]
and $\hat{x}$ commutes with $\hat{y} := y$.

If $|y|^2 \equiv 1 \pmod{4}$, then we imitate this proof
interchanging $x$ and $y$.

If $|x|^2, |y|^2 \equiv 3 \pmod{4}$,
then $x_0$, $z_x$, $c_2$, $c_3$, $y_0$, $z_y$, $d_2$, $d_3$ are odd, $c_1$, $d_1$ are even
and we can take 
\[
\hat{x} := x_0 + z_x(d_1 i + d_2 j + d_3 k) \in T_{p,l}
\]
and $\hat{y} := y$ as above.
\end{proof}

For an odd prime number $p$ we consider the set $n(X_p) = \{ n(x) : x \in X_p \}$. This
is a finite subset of $\mathbb{N}$ satisfying $\max n(X_p) \leq p-1$.

\begin{Corollary} \label{Cor26}
Let $p,l$ be two distinct odd prime numbers. 
There are commuting quaternions $x \in X_p$, $y \in X_l$, if and only if
$n(X_p) \cap n(X_l) \ne \emptyset$. 
\end{Corollary}

\begin{proof}
We combine Lemma~\ref{Lemma17} and Proposition~\ref{Prop25}
using that $(X_p \cup X_l) \subset T_{p,l}$ and
$(X_p \cup X_l) \cap \mathbb{Q} = \emptyset$.
\end{proof}

For $p < 200$,
we list all sets $n(X_p)$ in Table~\ref{Table2} in the Appendix. 
It shows for example that there are prime numbers
$p,l \equiv 7 \pmod{8}$ excluded from Proposition~\ref{Prop24} 
with commuting quaternions $x \in X_p$, $y \in X_l$.
For example take $p = 47$ and $l = 167$ 
(such that $l \equiv 23 \pmod{24}$ and $l \equiv 79 \pmod{88}$), 
where $n(X_{47}) \cap n(X_{167}) = \{46\}$,
and take commuting quaternions $x = 1+6i+j+3k$ and $y = 11+6i+j+3k$
such that $|x|^2 = 47$, $|y|^2 = 167$ and $n(x) = n(y) = 46$.

In Table~\ref{Table3} we also list $n(X_p)$ for all prime numbers $p < 1000$ 
satisfying $p \equiv 23 \pmod{24}$ and
$p \equiv 7, 39, 63, 79, 87 \pmod{88}$, that is for all prime
numbers $p \equiv 7 \pmod{8}$, $p < 1000$, excluded from Proposition~\ref{Prop24}(1).

\begin{Lemma} \label{Lemma27}
Let $p$ be an odd prime number and $m \in  n(X_p)$. Then there exist $r, s \in \mathbb{N}$ such
that $p = r^2 + m s^2$. 
\end{Lemma}

\begin{proof}
Let $x = x_0 + x_1 i + x_2 j + x_3 k \in X_p$ such that $n(x) = m$.
Since $n(x) = (x_1^2 + x_2^2 + x_3^2)/t^2$ 
for some $t \in \mathbb{N}$,
we have $p = x_0^2 + x_1^2 + x_2^2 + x_3^2 = x_0^2 + m t^2$, 
i.e.\ we can take $s = t$ and $r = |x_0|$.  
\end{proof}

The converse of Lemma~\ref{Lemma27} is true in some special cases
as seen in the proofs of Proposition~\ref{Prop23} and Proposition~\ref{Prop24}.
However it is not true in general, for example take $p = 7$, $m = 3$,
$r = 2$, $s = 1$ and observe that $m \notin n(X_7) = \{ 6 \}$.

\begin{Corollary} \label{Cor28}
Let $p,l$ be two distinct odd prime numbers. 
If there are commuting quaternions $x \in X_p$, $y \in X_l$, then
there exists an $m \in \mathbb{N}$
and $r_1, r_2, s_1, s_2 \in \mathbb{N}$
such that $p = r_1^2 + m s_1^2$ and $l = r_2^2 + m s_2^2$.
\end{Corollary}

\begin{proof}
Combine Corollary~\ref{Cor26} and Lemma~\ref{Lemma27}.
\end{proof}

For $p$ an odd prime number, let $n(X_p)_{min}$ 
be the smallest element in $n(X_p)$.
From what we have seen, it immediately follows that
$n(X_p)_{min} = 1$, if $p \equiv 1 \pmod{4}$, 
$n(X_p)_{min} = 2$, if $p \equiv 3 \pmod{8}$,
and 
$n(X_p)_{min} = 6$, if $p \equiv 7 \pmod{24}$.
In the remaining case $p \equiv 23 \pmod{24}$,
there seems to be no upper bound for $n(X_p)_{min}$.
We compute for example
$n(X_{23})_{min} = 14$, 
$n(X_{47})_{min} = 22$,
$n(X_{167})_{min} = 46$,
$n(X_{503})_{min} = 62$,
$n(X_{1223})_{min} = 134$,
$n(X_{1823})_{min} = 142$, 
$n(X_{1847})_{min} = 166$,
$n(X_{4703})_{min} = 214$,
$n(X_{8543})_{min} = 262$,
$n(X_{9743})_{min} = 334$.
This phenomenon makes it difficult to answer
Question~\ref{Question15} in full generality.

\newpage

\section{Appendix: some lists} \label{lists}
In the following list (Table~\ref{Table1}), we give some examples of pairs $x \in X_p$, $y \in X_l$
such that the index of $\langle \psi_{p,l}(x), \psi_{p,l}(y) \rangle$
in $\Gamma_{p,l}$ is finite (in particular Theorem~\ref{Thm11} can be applied). 
This index is given in the fifth column of the list.
We have used \textsf{GAP} (\cite{GAP}) for the computations.
We see that for $(p,l)=(3,5)$ there are only two possibilities for 
the index and abelianization. For the other pairs $(p,l)$, we 
have therefore only included some ``typical'' examples to keep the
list reasonably short. It also happened that for some 
non-commuting pairs $x$, $y$ (for example $x = 1+j+k$, $y = 1+6i+2k$,
or if $p$, $l$ are large) we were not able to compute
the index and the abelianization. In these cases we do
not know if the index is indeed infinite or finite
(but perhaps very large or difficult to compute).
Recall that for \emph{free} groups $\langle \psi_{p,l}(x), \psi_{p,l}(y) \rangle$,
the index would be infinite.
\setlongtables
\begin{longtable}{| r | r | c | c | r | c |} 
\hline
$p$ & $l$ & $x$     & $y$            & index & $\langle \psi_{p,l}(x), \psi_{p,l}(y) \rangle^{ab}$ \\ \hline \hline 
$3$ & $5$ & $1+j+k$ & $1+2i$         & $4$   & $\mathbb{Z}_{8} \times \mathbb{Z}_{16}$ \\ \hline
    &     & $1+j+k$ & $1+2j$         & $2$   & $\mathbb{Z}_{8} \times \mathbb{Z}_{8}$ \\ \hline
    &     & $1+j+k$ & $1+2k$         & $2$   & $\mathbb{Z}_{8} \times \mathbb{Z}_{8}$ \\ \hline
    &     & $1+j-k$ & $1+2i$         & $4$   & $\mathbb{Z}_{8} \times \mathbb{Z}_{16}$ \\ \hline
    &     & $1+j-k$ & $1+2j$         & $2$   & $\mathbb{Z}_{8} \times \mathbb{Z}_{8}$ \\ \hline
    &     & $1+j-k$ & $1+2k$         & $2$   & $\mathbb{Z}_{8} \times \mathbb{Z}_{8}$ \\ \hline
$3$ & $7$ & $1+j+k$ & $1+2i+j+k$     & $4$   & $\mathbb{Z}_{8} \times \mathbb{Z}_{16}$ \\ \hline
    &     & $1+j+k$ & $1+2i+j-k$     & $2$   & $\mathbb{Z}_{8} \times \mathbb{Z}_{8}$ \\ \hline
$3$ &$11$ & $1+j+k$ & $1+j+3k$       & $2$   & $\mathbb{Z}_{8} \times \mathbb{Z}_{8}$ \\ \hline
    &     & $1+j+k$ & $1+j-3k$       & $8$   & $\mathbb{Z}_{8} \times \mathbb{Z}_{32}$ \\ \hline
$3$ &$13$ & $1+j+k$ & $1+2i+2j+2k$   & $4$   & $\mathbb{Z}_{8} \times \mathbb{Z}_{16}$ \\ \hline
    &     & $1+j+k$ & $3+2i$         & $4$   & $\mathbb{Z}_{8} \times \mathbb{Z}_{16}$ \\ \hline
    &     & $1+j+k$ & $3+2j$         & $2$   & $\mathbb{Z}_{8} \times \mathbb{Z}_{8}$ \\ \hline
$3$ &$17$ & $1+j+k$ & $1+4i$         & $16$  & $\mathbb{Z}_{8} \times \mathbb{Z}_{64}$ \\ \hline
    &     & $1+j+k$ & $1+4j$         & $8$   & $\mathbb{Z}_{8} \times \mathbb{Z}_{32}$ \\ \hline
    &     & $1+j+k$ & $3+2i+2j$      & $2$   & $\mathbb{Z}_{8} \times \mathbb{Z}_{8}$ \\ \hline
    &     & $1+j+k$ & $3+2j-2k$      & $8$   & $\mathbb{Z}_{8} \times \mathbb{Z}_{32}$ \\ \hline
$3$ &$19$ & $1+j+k$ & $1+4i+j+k$     & $16$  & $\mathbb{Z}_{8} \times \mathbb{Z}_{64}$ \\ \hline
    &     & $1+j+k$ & $1+4i+j-k$     & $2$   & $\mathbb{Z}_{8} \times \mathbb{Z}_{8}$ \\ \hline
    &     & $1+j+k$ & $1+3j-3k$      & $2$   & $\mathbb{Z}_{8} \times \mathbb{Z}_{8}$ \\ \hline
    &     & $1+j+k$ & $3+j-3k$       & $8$   & $\mathbb{Z}_{8} \times \mathbb{Z}_{32}$ \\ \hline
$3$ &$23$ & $1+j+k$ & $1+2i+3j+3k$   & $4$   & $\mathbb{Z}_{8} \times \mathbb{Z}_{16}$ \\ \hline
    &     & $1+j+k$ & $1+2i+3j-3k$   & $48$  & $\mathbb{Z}_{8} \times \mathbb{Z}_{40}$ \\ \hline
    &     & $1+j+k$ & $3+2i+j+3k$    & $2$   & $\mathbb{Z}_{8} \times \mathbb{Z}_{8}$ \\ \hline
$3$ &$29$ & $1+j+k$ & $3+4i+2j$      & $2$   & $\mathbb{Z}_{8} \times \mathbb{Z}_{8}$ \\ \hline
    &     & $1+j+k$ & $3+2i+4j$      & $4$   & $\mathbb{Z}_{8} \times \mathbb{Z}_{16}$ \\ \hline
    &     & $1+j+k$ & $5+2i$         & $160$ & $\mathbb{Z}_{8} \times \mathbb{Z}_{16}$ \\ \hline
    &     & $1+j+k$ & $5+2j$         & $2$   & $\mathbb{Z}_{8} \times \mathbb{Z}_{8}$ \\ \hline
$3$ &$31$ & $1+j+k$ & $1+2i+j+5k$    & $4$   & $\mathbb{Z}_{8} \times \mathbb{Z}_{16}$ \\ \hline
    &     & $1+j+k$ & $1+2i+j-5k$    & $48$  & $\mathbb{Z}_{8} \times \mathbb{Z}_{40}$ \\ \hline
    &     & $1+j+k$ & $5+2i+j-k$     & $80$  & $\mathbb{Z}_{8} \times \mathbb{Z}_{8}$ \\ \hline
    &     & $1+j+k$ & $3+2i+3j+3k$   & $4$   & $\mathbb{Z}_{8} \times \mathbb{Z}_{16}$ \\ \hline
    &     & $1+j+k$ & $3+2i+3j-3k$   & $48$  & $\mathbb{Z}_{8} \times \mathbb{Z}_{40}$ \\ \hline
$3$ &$37$ & $1+j+k$ & $1+6i$         & $4$   & $\mathbb{Z}_{8} \times \mathbb{Z}_{16}$ \\ \hline
    &     & $1+j+k$ & $1+6j$         & $22$  & $\mathbb{Z}_{8} \times \mathbb{Z}_{8}$ \\ \hline
    &     & $1+j+k$ & $1+2i+4j+4k$   & $4$   & $\mathbb{Z}_{8} \times \mathbb{Z}_{16}$ \\ \hline
    &     & $1+j+k$ & $1+4i+2j+4k$   & $22$  & $\mathbb{Z}_{8} \times \mathbb{Z}_{8}$ \\ \hline
    &     & $1+j+k$ & $1+4i+2j-4k$   & $144$ & $\mathbb{Z}_{8} \times \mathbb{Z}_{32}$ \\ \hline
    &     & $1+j+k$ & $5+2i+2j+2k$   & $4$   & $\mathbb{Z}_{8} \times \mathbb{Z}_{16}$ \\ \hline
    &     & $1+j+k$ & $5+2i+2j-2k$   & $160$ & $\mathbb{Z}_{8} \times \mathbb{Z}_{16}$ \\ \hline
$3$ &$41$ & $1+j+k$ & $1+2j+6k$      & $8$   & $\mathbb{Z}_{8} \times \mathbb{Z}_{32}$ \\ \hline
    &     & $1+j+k$ & $1+2j-6k$      & $32$  & $\mathbb{Z}_{8} \times \mathbb{Z}_{128}$ \\ \hline
    &     & $1+j+k$ & $1+2i+6k$      & $48$  & $\mathbb{Z}_{8} \times \mathbb{Z}_{40}$ \\ \hline
    &     & $1+j+k$ & $3+4i+4j$      & $8$   & $\mathbb{Z}_{8} \times \mathbb{Z}_{32}$ \\ \hline
    &     & $1+j+k$ & $3+4i-4k$      & $32$  & $\mathbb{Z}_{8} \times \mathbb{Z}_{128}$ \\ \hline
    &     & $1+j+k$ & $5+4j$         & $8$   & $\mathbb{Z}_{8} \times \mathbb{Z}_{32}$ \\ \hline
$3$ &$43$ & $1+j+k$ & $3+3j-5k$      & $32$  & $\mathbb{Z}_{8} \times \mathbb{Z}_{128}$ \\ \hline
    &     & $1+j+k$ & $5+3j-3k$      & $80$  & $\mathbb{Z}_{8} \times \mathbb{Z}_{8}$ \\ \hline
    &     & $1+j+k$ & $1+4i+j+5k$    & $8$   & $\mathbb{Z}_{8} \times \mathbb{Z}_{32}$ \\ \hline
    &     & $1+j+k$ & $1+4i+j-5k$    & $144$ & $\mathbb{Z}_{8} \times \mathbb{Z}_{32}$ \\ \hline
    &     & $1+j+k$ & $5+4i+j+k$     & $16$  & $\mathbb{Z}_{8} \times \mathbb{Z}_{64}$ \\ \hline
    &     & $1+j+k$ & $5+4i+j-k$     & $80$  & $\mathbb{Z}_{8} \times \mathbb{Z}_{8}$ \\ \hline
    &     & $1+j+k$ & $3+4i+3j+3k$   & $16$  & $\mathbb{Z}_{8} \times \mathbb{Z}_{64}$ \\ \hline
    &     & $1+j+k$ & $3+4i+3j-3k$   & $144$ & $\mathbb{Z}_{8} \times \mathbb{Z}_{32}$ \\ \hline
$3$ &$47$ & $1+j+k$ & $3+2i+3j-5k$   & $4$   & $\mathbb{Z}_{8} \times \mathbb{Z}_{16}$ \\ \hline
    &     & $1+j+k$ & $5+2i+3j-3k$   & $1920$& $\mathbb{Z}_{8} \times \mathbb{Z}_{40}$ \\ \hline
    &     & $1+j+k$ & $1+6i+j+3k$    & $80$  & $\mathbb{Z}_{8} \times \mathbb{Z}_{72}$ \\ \hline
    &     & $1+j+k$ & $1+6i+j-3k$    & $96$  & $\mathbb{Z}_{8} \times \mathbb{Z}_{80}$ \\ \hline
    &     & $1+j+k$ & $3+6i+j+k$     & $4$   & $\mathbb{Z}_{8} \times \mathbb{Z}_{16}$ \\ \hline
$3$ &$53$ & $1+j+k$ & $1+4j+6k$      & $22$  & $\mathbb{Z}_{8} \times \mathbb{Z}_{8}$ \\ \hline
    &     & $1+j+k$ & $3+2i+2j+6k$   & $4$   & $\mathbb{Z}_{8} \times \mathbb{Z}_{16}$ \\ \hline
    &     & $1+j+k$ & $7+2k$         & $22$  & $\mathbb{Z}_{8} \times \mathbb{Z}_{8}$ \\ \hline
$3$ &$59$ & $1+j+k$ & $5+4i+3j+3k$   & $16$  & $\mathbb{Z}_{8} \times \mathbb{Z}_{64}$ \\ \hline
$3$ &$61$ & $1+j+k$ & $3+4j+6k$      & $2$   & $\mathbb{Z}_{8} \times \mathbb{Z}_{8}$ \\ \hline
    &     & $1+j+k$ & $3+4j-6k$      & $80$  & $\mathbb{Z}_{24} \times \mathbb{Z}_{40}$ \\ \hline
    &     & $1+j+k$ & $5+6k$         & $2$   & $\mathbb{Z}_{8} \times \mathbb{Z}_{8}$ \\ \hline
    &     & $1+j+k$ & $5+4i+2j+4k$   & $2$   & $\mathbb{Z}_{8} \times \mathbb{Z}_{8}$ \\ \hline
$3$ &$67$ & $1+j+k$ & $1+4i+j+7k$    & $144$ & $\mathbb{Z}_{8} \times \mathbb{Z}_{32}$ \\ \hline
    &     & $1+j+k$ & $1+4i+j-7k$    & $16$  & $\mathbb{Z}_{8} \times \mathbb{Z}_{64}$ \\ \hline
    &     & $1+j+k$ & $5+4i+j-5k$    & $144$ & $\mathbb{Z}_{8} \times \mathbb{Z}_{32}$ \\ \hline
$3$ &$71$ & $1+j+k$ & $3+6i+j+5k$    & $96$  & $\mathbb{Z}_{8} \times \mathbb{Z}_{80}$ \\ \hline
$3$ &$73$ & $1+j+k$ & $1+6j-6k$      & $8$   & $\mathbb{Z}_{8} \times \mathbb{Z}_{32}$ \\ \hline
    &     & $1+j+k$ & $1+6i+6k$      & $2$   & $\mathbb{Z}_{8} \times \mathbb{Z}_{8}$ \\ \hline
    &     & $1+j+k$ & $1+2i+2j+8k$   & $48$  & $\mathbb{Z}_{8} \times \mathbb{Z}_{40}$ \\ \hline
    &     & $1+j+k$ & $1+2i+2j-8k$   & $2$   & $\mathbb{Z}_{8} \times \mathbb{Z}_{8}$ \\ \hline
    &     & $1+j+k$ & $1+8i+2j+2k$   & $64$  & $\mathbb{Z}_{8} \times \mathbb{Z}_{256}$ \\ \hline
    &     & $1+j+k$ & $3+8i$         & $64$  & $\mathbb{Z}_{8} \times \mathbb{Z}_{256}$ \\ \hline
    &     & $1+j+k$ & $3+8j$         & $32$  & $\mathbb{Z}_{8} \times \mathbb{Z}_{128}$ \\ \hline
    &     & $1+j+k$ & $5+4i+4j+4k$   & $16$  & $\mathbb{Z}_{8} \times \mathbb{Z}_{64}$ \\ \hline
    &     & $1+j+k$ & $5+4i+4j-4k$   & $640$ & $\mathbb{Z}_{8} \times \mathbb{Z}_{64}$ \\ \hline
    &     & $1+j+k$ & $7+2i+2j+4k$   & $2$   & $\mathbb{Z}_{8} \times \mathbb{Z}_{8}$ \\ \hline
    &     & $1+j+k$ & $7+2i+2j-4k$   & $48$  & $\mathbb{Z}_{8} \times \mathbb{Z}_{40}$ \\ \hline
    &     & $1+j+k$ & $7+4i+2j+2k$   & $16$  & $\mathbb{Z}_{8} \times \mathbb{Z}_{64}$ \\ \hline
    &     & $1+j+k$ & $7+4i+2j-2k$   & $672$ & $\mathbb{Z}_{8} \times \mathbb{Z}_{64}$ \\ \hline
$3$ &$79$ & $1+j+k$ & $3+6i+3j+5k$   & $80$  & $\mathbb{Z}_{8} \times \mathbb{Z}_{72}$ \\ \hline
    &     & $1+j+k$ & $5+2i+j-7k$    & $4$   & $\mathbb{Z}_{8} \times \mathbb{Z}_{16}$ \\ \hline
    &     & $1+j+k$ & $5+6i+3j+3k$   & $4$   & $\mathbb{Z}_{8} \times \mathbb{Z}_{16}$ \\ \hline
$3$ &$83$ & $1+j+k$ & $1+8i+3j+3k$   & $64$  & $\mathbb{Z}_{8} \times \mathbb{Z}_{256}$ \\ \hline
$3$ &$89$ & $1+j+k$ & $3+4i+8k$      & $16$  & $\mathbb{Z}_{8} \times \mathbb{Z}_{64}$ \\ \hline\hline
$5$ & $7$ & $1+2i$  & $1+2i+j+k$     & $4$   & $\mathbb{Z}_{8} \times \mathbb{Z}_{16}$ \\ \hline
    &     & $1+2j$  & $1+2i+j+k$     & $2$   & $\mathbb{Z}_{8} \times \mathbb{Z}_{8}$ \\ \hline
$5$ & $11$& $1+2i$  & $1+j+3k$       & $4$   & $\mathbb{Z}_{8} \times \mathbb{Z}_{16}$ \\ \hline
    &     & $1+2i$  & $3+j+k$        & $48$  & $\mathbb{Z}_{16} \times \mathbb{Z}_{16}$ \\ \hline
    &     & $1+2j$  & $1+j+3k$       & $48$  & $\mathbb{Z}_{8} \times \mathbb{Z}_{24}$ \\ \hline
    &     & $1+2j$  & $1+3j+k$       & $2$   & $\mathbb{Z}_{8} \times \mathbb{Z}_{8}$ \\ \hline
$5$ & $13$& $1+2i$  & $1+2i+2j+2k$   & $16$  & $\mathbb{Z}_{16} \times \mathbb{Z}_{32}$ \\ \hline
    &     & $1+2i$  & $3+2j$         & $96$  & $\mathbb{Z}_{16} \times \mathbb{Z}_{32}$ \\ \hline
$5$ & $17$& $1+2i$  & $1+4j$         & $32$  & $\mathbb{Z}_{16} \times \mathbb{Z}_{64}$ \\ \hline
    &     & $1+2i$  & $3+2i+2j$      & $8$   & $\mathbb{Z}_{16} \times \mathbb{Z}_{16}$ \\ \hline
    &     & $1+2i$  & $3+2j+2k$      & $192$ & $\mathbb{Z}_{16} \times \mathbb{Z}_{64}$ \\ \hline
$5$ & $19$& $1+2i$  & $1+4i+j+k$     & $4$   & $\mathbb{Z}_{8} \times \mathbb{Z}_{16}$ \\ \hline
    &     & $1+2i$  & $1+3j+3k$      & $96$  & $\mathbb{Z}_{8} \times \mathbb{Z}_{48}$ \\ \hline
    &     & $1+2i$  & $3+j+3k$       & $48$  & $\mathbb{Z}_{16} \times \mathbb{Z}_{16}$ \\ \hline
    &     & $1+2j$  & $1+4i+j+k$     & $144$ & $\mathbb{Z}_{8} \times \mathbb{Z}_{32}$ \\ \hline
    &     & $1+2j$  & $1+3j+3k$      & $48$  & $\mathbb{Z}_{8} \times \mathbb{Z}_{24}$ \\ \hline
    &     & $1+2j$  & $3+3j+k$       & $24$  & $\mathbb{Z}_{8} \times \mathbb{Z}_{16}$ \\ \hline
$5$ & $23$& $1+2i$  & $1+2i+3j+3k$   & $96$  & $\mathbb{Z}_{8} \times \mathbb{Z}_{48}$ \\ \hline
    &     & $1+2i$  & $3+2i+j+3k$    & $4$   & $\mathbb{Z}_{8} \times \mathbb{Z}_{16}$ \\ \hline
    &     & $1+2j$  & $1+2i+3j+3k$   & $112$ & $\mathbb{Z}_{8} \times \mathbb{Z}_{24}$ \\ \hline
    &     & $1+2j$  & $3+2i+3j+k$    & $24$  & $\mathbb{Z}_{8} \times \mathbb{Z}_{16}$ \\ \hline
$5$ & $29$& $1+2i$  & $3+4i+2j$      & $8$   & $\mathbb{Z}_{16} \times \mathbb{Z}_{16}$ \\ \hline
    &     & $1+2i$  & $3+2i+4j$      & $32$  & $\mathbb{Z}_{16} \times \mathbb{Z}_{64}$ \\ \hline
    &     & $1+2i$  & $3+2j+4k$      & $96$  & $\mathbb{Z}_{16} \times \mathbb{Z}_{32}$ \\ \hline
    &     & $1+2i$  & $5+2j$         & $8$   & $\mathbb{Z}_{16} \times \mathbb{Z}_{16}$ \\ \hline
$5$ & $31$& $1+2i$  & $1+2i+j+5k$    & $224$ & $\mathbb{Z}_{8} \times \mathbb{Z}_{48}$ \\ \hline
    &     & $1+2j$  & $1+2i+j+5k$    & $240$ & $\mathbb{Z}_{8} \times \mathbb{Z}_{56}$ \\ \hline
    &     & $1+2j$  & $3+2i+3j+3k$   & $1344$& $\mathbb{Z}_{8} \times \mathbb{Z}_{48}$ \\ \hline
$5$ & $37$& $1+2i$  & $1+4i+2j+4k$   & $8$   & $\mathbb{Z}_{16} \times \mathbb{Z}_{16}$ \\ \hline
$5$ & $41$& $1+2i$  & $1+2j+6k$      & $16$  & $\mathbb{Z}_{16} \times \mathbb{Z}_{32}$ \\ \hline
    &     & $1+2i$  & $1+2i+6k$      & $192$ & $\mathbb{Z}_{16} \times \mathbb{Z}_{48}$ \\ \hline
    &     & $1+2i$  & $3+4j+4k$      & $768$ & $\mathbb{Z}_{16} \times \mathbb{Z}_{256}$ \\ \hline
    &     & $1+2i$  & $3+4i+4k$      & $32$  & $\mathbb{Z}_{16} \times \mathbb{Z}_{64}$ \\ \hline
    &     & $1+2i$  & $5+4j$         & $32$  & $\mathbb{Z}_{16} \times \mathbb{Z}_{64}$ \\ \hline 
$5$ & $43$& $1+2j$  & $3+3j+5k$      & $24$  & $\mathbb{Z}_{8} \times \mathbb{Z}_{16}$ \\ \hline
    &     & $1+2j$  & $3+5j+3k$      & $48$  & $\mathbb{Z}_{8} \times \mathbb{Z}_{24}$ \\ \hline
    &     & $1+2j$  & $3+4i+3j+3k$   & $24$  & $\mathbb{Z}_{8} \times \mathbb{Z}_{16}$ \\ \hline
$5$ & $47$& $1+2i$  & $3+2i+5j-3k$   & $288$ & $\mathbb{Z}_{16} \times \mathbb{Z}_{32}$ \\ \hline
    &     & $1+2j$  & $3+2i+5j+3k$   & $112$ & $\mathbb{Z}_{8} \times \mathbb{Z}_{24}$ \\ \hline\hline
$7$ & $11$& $1+2i+j+k$  & $1+j+3k$   & $2$   & $\mathbb{Z}_{8} \times \mathbb{Z}_{8}$ \\ \hline
    &     & $1+2i+j+k$  & $1+j-3k$   & $4$   & $\mathbb{Z}_{8} \times \mathbb{Z}_{16}$ \\ \hline
    &     & $1+2i+j+k$  & $3+j+k$    & $20$  & $\mathbb{Z}_{8} \times \mathbb{Z}_{16}$ \\ \hline
    &     & $1+2i+j+k$  & $3+j-k$    & $48$  & $\mathbb{Z}_{8} \times \mathbb{Z}_{24}$ \\ \hline
$7$ & $13$& $1+2i+j+k$  & $1+2i+2j+2k$ & $4$ & $\mathbb{Z}_{8} \times \mathbb{Z}_{48}$ \\ \hline
    &     & $1+2i+j+k$  & $1-2i+2j+2k$ & $32$& $\mathbb{Z}_{24} \times \mathbb{Z}_{48}$ \\ \hline
    &     & $1+2i+j+k$  & $3+2i$     & $4$   & $\mathbb{Z}_{8} \times \mathbb{Z}_{48}$ \\ \hline
    &     & $1+2i+j+k$  & $3+2j$     & $48$  & $\mathbb{Z}_{8} \times \mathbb{Z}_{24}$ \\ \hline
$7$ & $17$& $1+2i+j+k$  & $1+4i$     & $16$  & $\mathbb{Z}_{8} \times \mathbb{Z}_{64}$ \\ \hline
    &     & $1+2i+j+k$  & $1+4j$     & $192$ & $\mathbb{Z}_{8} \times \mathbb{Z}_{32}$ \\ \hline
    &     & $1+2i+j+k$  & $3+2i+2j$  & $48$  & $\mathbb{Z}_{8} \times \mathbb{Z}_{24}$ \\ \hline
    &     & $1+2i+j+k$  & $3+2i-2j$  & $2736$& $\mathbb{Z}_{8} \times \mathbb{Z}_{40}$ \\ \hline
    &     & $1+2i+j+k$  & $3+2j+2k$  & $16$  & $\mathbb{Z}_{8} \times \mathbb{Z}_{64}$ \\ \hline
    &     & $1+2i+j+k$  & $3+2j-2k$  & $192$&  $\mathbb{Z}_{8} \times \mathbb{Z}_{96}$ \\ \hline
$7$ & $19$& $1+2i+j+k$  & $1+4i+j+k$ & $4$   & $\mathbb{Z}_{8} \times \mathbb{Z}_{48}$ \\ \hline
    &     & $1+2i+j+k$  & $1+4i+j-k$ & $48$  & $\mathbb{Z}_{24} \times \mathbb{Z}_{40}$ \\ \hline
    &     & $1+2i+j+k$  & $1+4i-j-k$ & $160$ & $\mathbb{Z}_{24} \times \mathbb{Z}_{48}$ \\ \hline
    &     & $1+2i+j+k$  & $1+3j+3k$  & $160$ & $\mathbb{Z}_{24} \times \mathbb{Z}_{48}$ \\ \hline
    &     & $1+2i+j+k$  & $3+j+3k$   & $48$  & $\mathbb{Z}_{24} \times \mathbb{Z}_{40}$ \\ \hline
    &     & $1+2i+j+k$  & $3+j-3k$   & $4$   & $\mathbb{Z}_{8} \times \mathbb{Z}_{48}$ \\ \hline
$7$ & $23$& $1+2i+j+k$  & $3+2i+3j-k$& $192$ & $\mathbb{Z}_{8} \times \mathbb{Z}_{96}$ \\ \hline
$7$ & $29$& $1+2i+j+k$  & $3+4i+2j$  & $48$  & $\mathbb{Z}_{8} \times \mathbb{Z}_{8}$ \\ \hline
    &     & $1+2i+j+k$  & $3+2j-4k$  & $2$   & $\mathbb{Z}_{8} \times \mathbb{Z}_{8}$ \\ \hline
    &     & $1+2i+j+k$  & $5+2j$     & $48$  & $\mathbb{Z}_{8} \times \mathbb{Z}_{8}$ \\ \hline
$7$ & $31$& $1+2i+j+k$  & $1+2i-j+5k$& $240$ & $\mathbb{Z}_{24} \times \mathbb{Z}_{56}$ \\ \hline
    &     & $1+2i+j+k$  & $3+2i+3j-3k$& $240$& $\mathbb{Z}_{24} \times \mathbb{Z}_{56}$ \\ \hline\hline
$11$& $13$&$1+j+3k$     & $1+2i+2j-2k$ & $4$ & $\mathbb{Z}_{8} \times \mathbb{Z}_{16}$ \\ \hline
    &     &$1+j+3k$     & $3+2k$     & $72$  & $\mathbb{Z}_{8} \times \mathbb{Z}_{16}$ \\ \hline
    &     &$3+j+k$      & $1+2i+2j-2k$ &$288$& $\mathbb{Z}_{8} \times \mathbb{Z}_{48}$ \\ \hline
    &     &$3+j+k$      & $3+2j$     & $72$  & $\mathbb{Z}_{8} \times \mathbb{Z}_{16}$ \\ \hline
    &     &$3+j-k$      & $1+2i-2j+2k$ & $4$ & $\mathbb{Z}_{8} \times \mathbb{Z}_{16}$ \\ \hline
$11$& $17$&$1+j+3k$     & $1+4j$     & $192$ & $\mathbb{Z}_{8} \times \mathbb{Z}_{96}$ \\ \hline
    &     &$1+j+3k$     & $1+4k$     & $8$   & $\mathbb{Z}_{8} \times \mathbb{Z}_{32}$ \\ \hline
    &     &$1+j+3k$     & $3+2i+2k$  & $72$  & $\mathbb{Z}_{8} \times \mathbb{Z}_{16}$ \\ \hline
    &     &$1+j+3k$     & $3+2j+2k$  & $8$   & $\mathbb{Z}_{8} \times \mathbb{Z}_{32}$ \\ \hline
    &     &$3+j+k$      & $3+2i+2j$  & $144$ & $\mathbb{Z}_{8} \times \mathbb{Z}_{24}$ \\ \hline
$11$& $19$&$3+j+k$      & $1+4i+j+k$ & $80$  & $\mathbb{Z}_{8} \times \mathbb{Z}_{192}$ \\ \hline
    &     &$3+j+k$      & $3+j-3k$   & $96$  & $\mathbb{Z}_{16} \times \mathbb{Z}_{32}$ \\ \hline
$11$& $23$&$1+j+3k$     & $2+i+3j+3k$& $80$  & $\mathbb{Z}_{8} \times \mathbb{Z}_{72}$ \\ \hline
    &     &$1+j+3k$     & $2+i+3j-3k$& $56$  & $\mathbb{Z}_{8} \times \mathbb{Z}_{16}$ \\ \hline
    &     &$3+j+k$      & $2+i+3j-3k$& $336$ & $\mathbb{Z}_{8} \times \mathbb{Z}_{16}$ \\ \hline
    &     &$3+j+k$      & $2+3i+j-3k$& $96$  & $\mathbb{Z}_{8} \times \mathbb{Z}_{48}$ \\ \hline
$11$& $29$&$1+j+3k$     & $3+2j+4k$  & $2$   & $\mathbb{Z}_{8} \times \mathbb{Z}_{8}$ \\ \hline
    &     &$1+j+3k$     & $5+2k$     & $2$   & $\mathbb{Z}_{8} \times \mathbb{Z}_{8}$ \\ \hline
$11$& $31$&$3+j+k$      & $1+2i+j-5k$& $2$   & $\mathbb{Z}_{8} \times \mathbb{Z}_{8}$ \\ \hline
$11$& $37$&$1+j+3k$     & $1+2i+4j-4k$& $56$ & $\mathbb{Z}_{8} \times \mathbb{Z}_{16}$ \\ \hline
    &     &$1+j+3k$     & $1+4i+2j+4k$& $336$& $\mathbb{Z}_{8} \times \mathbb{Z}_{80}$ \\ \hline
$11$& $41$&$3+j+k$      & $1+2j+6k$&   $8$   & $\mathbb{Z}_{8} \times \mathbb{Z}_{32}$ \\ \hline
    &     &$3+j+k$      & $5+4j$&      $8$   & $\mathbb{Z}_{8} \times \mathbb{Z}_{32}$ \\ \hline \hline
$13$& $17$&$3+2i$       & $3+2i+2j$&   $288$ & $\mathbb{Z}_{16} \times \mathbb{Z}_{32}$ \\ \hline
$13$& $19$&$1+2i+2j+2k$ & $1+4i+j-k$&  $4$   & $\mathbb{Z}_{8} \times \mathbb{Z}_{48}$ \\ \hline
    &     &$1+2i+2j+2k$ & $3+j-3k$&    $4$   & $\mathbb{Z}_{8} \times \mathbb{Z}_{48}$ \\ \hline
$13$& $23$&$3+2i$       & $3+2i+j+3k$& $1152$& $\mathbb{Z}_{8} \times \mathbb{Z}_{32}$ \\ \hline
    &     &$3+2j$       & $1+2i+3j+3k$&$24$  & $\mathbb{Z}_{8} \times \mathbb{Z}_{16}$ \\ \hline
$13$& $29$&$3+2i$       & $3+2i+4j$ &  $384$ & $\mathbb{Z}_{32} \times \mathbb{Z}_{64}$ \\ \hline
$13$& $31$&$1+2i+2j+2k$ & $3+2i+3j-3k$ & $4$ & $\mathbb{Z}_{8} \times \mathbb{Z}_{48}$ \\ \hline \hline
$17$& $19$&$1+4i$       & $1+4i+j+k$& $1920$ & $\mathbb{Z}_{8} \times \mathbb{Z}_{64}$ \\ \hline
    &     &$3+2i+2j$    & $3+j+3k$  & $72$   & $\mathbb{Z}_{8} \times \mathbb{Z}_{16}$ \\ \hline
    &     &$3+2j+2k$    & $1+4i+j-k$ & $576$ & $\mathbb{Z}_{8} \times \mathbb{Z}_{288}$ \\ \hline
    &     &$3+2j+2k$    & $1+3j-3k$  & $576$ & $\mathbb{Z}_{8} \times \mathbb{Z}_{288}$ \\ \hline
    &     &$3+2j+2k$    & $3+j+3k$   & $96$  & $\mathbb{Z}_{8} \times \mathbb{Z}_{64}$ \\ \hline \hline 
$19$& $23$&$3+j+3k$     & $1+2i+3j-3k$& $48$ & $\mathbb{Z}_{8} \times \mathbb{Z}_{32}$ \\ \hline
    &     &$3+j+3k$     & $3+2i+j-3k$&  $24$ & $\mathbb{Z}_{8} \times \mathbb{Z}_{16}$ \\ \hline
\multicolumn{6}{c}{\phantom{a}} \\
\caption{Some examples where $[\Gamma_{p,l} : \langle \psi_{p,l}(x), \psi_{p,l}(y) \rangle] < \infty$} \label{Table1}
\end{longtable}

\newpage

By Corollary~\ref{Cor26}, the following table can be used to
answer Question~\ref{Question15} for all pairs $p,l < 200$.
\begin{table}[ht]
\begin{tabular}{r|l||r|l}
$p$ & $n(X_p)$ & $p$ & $n(X_p)$ \\ \hline
$3$ &  $\{  2 \}$ & $101$ &   $\{  1, 5, 13, 19, 25 \}$\\
$5$ &  $\{  1 \}$ & $103$ &   $\{  6, 22, 54, 78, 94, 102 \}$\\
$7$ &  $\{  6 \}$ & $107$ &   $\{  2, 26, 58, 82, 98, 106 \}$\\
$11$ &  $\{  2, 10 \}$ & $109$ &   $\{  1, 3, 21, 25, 27 \}$\\
$13$ &  $\{  1, 3 \}$ & $113$ &   $\{  1, 2, 22, 26 \}$\\
$17$ &  $\{  1, 2 \}$ & $127$ &   $\{  6, 14, 46, 78, 102, 118, 126 \}$\\
$19$ &  $\{  2, 10, 18 \}$ & $131$ &  $\{  2, 10, 50, 82, 106, 122, 130 \}$\\
$23$ &  $\{  14, 22 \}$ & $137$ &  $\{  1, 2, 14, 22, 34 \}$\\
$29$ &  $\{  1, 5 \}$ & $139$ &  $\{  2, 10, 18, 58, 90, 114, 130, 138 \}$\\
$31$ &  $\{  6, 22, 30 \}$ & $149$ &  $\{  1, 17, 25, 35, 37 \}$\\
$37$ &  $\{  1, 3, 9 \}$ & $151$ &  $\{  6, 14, 30, 70, 102, 126, 142, 150 \}$\\
$41$ &  $\{  1, 2, 10 \}$ & $157$ &  $\{  1, 3, 9, 19, 27, 33, 37 \}$\\
$43$ &  $\{  2, 18, 34, 42 \}$ & $163$ &  $\{  2, 18, 42, 82, 114, 138, 154, 162 \}$\\
$47$ &   $\{  22, 38, 46 \}$ & $167$ &  $\{  46, 86, 118, 142, 158, 166 \}$\\
$53$ &   $\{  1, 11, 13 \}$ & $173$ &  $\{  1, 13, 37, 41, 43 \}$\\
$59$ &   $\{  2, 10, 34, 50, 58 \}$ & $179$ &  $\{  2, 10, 58, 98, 130, 154, 170, 178 \}$\\
$61$ &  $\{  1, 3, 9, 13 \}$ & $181$ &  $\{  1, 3, 5, 25, 33, 43, 45 \}$\\
$67$ &  $\{  2, 18, 42, 58, 66 \}$ & $191$ &  $\{  22, 70, 110, 142, 166, 182, 190 \}$\\
$71$ &  $\{  22, 46, 62, 70 \}$ & $193$ &  $\{  1, 2, 3, 6, 9, 18, 42, 46 \}$\\
$73$ &  $\{  1, 2, 3, 6, 18 \}$ & $197$ &  $\{  1, 19, 29, 37, 43, 49 \}$\\
$79$ &   $\{  6, 30, 54, 70, 78 \}$ & $199$ &  $\{  6, 22, 30, 78, 118, 150, 174, 190, 198 \}$\\
$83$ &   $\{  2, 34, 58, 74, 82 \}$ & &\\
$89$ &   $\{  1, 2, 5, 10, 22 \}$ & &\\
$97$ &   $\{  1, 2, 3, 6, 18, 22 \}$ & &\\
\multicolumn{4}{c}{\phantom{a}} \\
\end{tabular}
\caption{$n(X_p)$ for prime numbers $2 < p < 200$} \label{Table2}
\end{table}

\newpage
Table~\ref{Table2b} contains the same information as Table~\ref{Table2},
but is sorted by $p \pmod{8}$. It is a good illustration of Lemma~\ref{Lemma19}. 
\begin{table}[ht]
\begin{tabular}{r|l||r|l}
$p$ & $n(X_p)$ & $p$ & $n(X_p)$ \\ \hline
$3$ &  $\{  2 \}$ & $7$ &  $\{  6 \}$ \\
$11$ &  $\{  2, 10 \}$ & $23$ &  $\{  14, 22 \}$ \\
$19$ &  $\{  2, 10, 18 \}$ & $31$ &  $\{  6, 22, 30 \}$ \\
$43$ &  $\{  2, 18, 34, 42 \}$ & $47$ &   $\{  22, 38, 46 \}$ \\
$59$ &   $\{  2, 10, 34, 50, 58 \}$ & $71$ &  $\{  22, 46, 62, 70 \}$ \\
$67$ &  $\{  2, 18, 42, 58, 66 \}$ & $79$ &   $\{  6, 30, 54, 70, 78 \}$ \\
$83$ &   $\{  2, 34, 58, 74, 82 \}$ & $103$ &   $\{  6, 22, 54, 78, 94, 102 \}$ \\ 
$107$ &   $\{  2, 26, 58, 82, 98, 106 \}$ & $127$ &   $\{  6, 14, 46, 78, 102, 118, 126 \}$ \\
$131$ &  $\{  2, 10, 50, 82, 106, 122, 130 \}$ & $151$ &  $\{  6, 14, 30, 70, 102, 126, 142, 150 \}$ \\
$139$ &  $\{  2, 10, 18, 58, 90, 114, 130, 138 \}$ & $167$ &  $\{  46, 86, 118, 142, 158, 166 \}$ \\
$163$ &  $\{  2, 18, 42, 82, 114, 138, 154, 162 \}$ & $191$ &  $\{  22, 70, 110, 142, 166, 182, 190 \}$ \\
$179$ &  $\{  2, 10, 58, 98, 130, 154, 170, 178 \}$ & $199$ &  $\{  6, 22, 30, 78, 118, 150, 174, 190, 198 \}$ \\ \hline
$5$ &  $\{  1 \}$ & $17$ &  $\{  1, 2 \}$ \\
$13$ &  $\{  1, 3 \}$ & $41$ &  $\{  1, 2, 10 \}$ \\
$29$ &  $\{  1, 5 \}$ & $73$ &  $\{  1, 2, 3, 6, 18 \}$ \\
$37$ &  $\{  1, 3, 9 \}$ & $89$ &   $\{  1, 2, 5, 10, 22 \}$ \\
$53$ &   $\{  1, 11, 13 \}$ & $97$ &   $\{  1, 2, 3, 6, 18, 22 \}$ \\
$61$ &  $\{  1, 3, 9, 13 \}$ & $113$ &   $\{  1, 2, 22, 26 \}$ \\
$101$ &   $\{  1, 5, 13, 19, 25 \}$ & $137$ &  $\{  1, 2, 14, 22, 34 \}$ \\
$109$ &   $\{  1, 3, 21, 25, 27 \}$ & $193$ &  $\{  1, 2, 3, 6, 9, 18, 42, 46 \}$ \\
$149$ &  $\{  1, 17, 25, 35, 37 \}$ & & \\
$157$ &  $\{  1, 3, 9, 19, 27, 33, 37 \}$ & & \\
$173$ &  $\{  1, 13, 37, 41, 43 \}$ & & \\
$181$ &  $\{  1, 3, 5, 25, 33, 43, 45 \}$ & & \\
$197$ &  $\{  1, 19, 29, 37, 43, 49 \}$ & & \\
\multicolumn{4}{c}{\phantom{a}} \\
\end{tabular}
\caption{$n(X_p)$ for prime numbers $2 < p < 200$, sorted by $p \pmod{8}$} \label{Table2b}
\end{table}

\newpage
In Table~\ref{Table3} we list $n(X_p)$ for all prime numbers $p < 1000$ 
satisfying $p \equiv 23 \pmod{24}$ and
$p \equiv 7, 39, 63, 79, 87 \pmod{88}$.
\begin{table}[ht]
\begin{tabular}{r|l}
$p$ & $n(X_p)$ \\ \hline
$167$ & $\{  46, 86, 118, 142, 158, 166 \}$\\
$239$ & $\{ 14, 70, 118, 158, 190, 214, 230, 238 \}$\\
$263$ & $\{ 38, 94, 142, 182, 214, 238, 254, 262 \}$\\
$359$ & $\{ 14, 70, 134, 190, 238, 278, 310, 334, 350, 358 \}$\\
$431$ & $\{ 14, 70, 142, 206, 262, 310, 350, 382, 406, 422, 430 \}$\\
$479$ & $\{ 38, 118, 190, 254, 310, 358, 398, 430, 454, 470, 478 \}$\\
$503$ & $\{ 62, 142, 214, 278, 334, 382, 422, 454, 478, 494, 502 \}$\\
$743$ & $\{ 14, 118, 214, 302, 382, 454, 518, 574, 622, 662, 694, 718, 734, 742 \}$\\
$887$ & $\{ 46, 158, 262, 358, 446, 526, 598, 662, 718, 766, 806, 838, 862, 878, 886 \}$\\
\multicolumn{2}{c}{\phantom{a}} \\
\end{tabular}
\caption{$n(X_p)$ for some prime numbers $p \equiv 23 \pmod{24}$} \label{Table3}
\end{table}


\newpage

\begin{thebibliography}{99}
\bibitem{BMII}
  Burger, Marc; Mozes, Shahar,
  \emph{Lattices in product of trees},
  Inst. Hautes \'Etudes Sci. Publ. Math. No. \textbf{92}(2000),
  151--194(2001).
\bibitem{Cox}
   Cox, David A., 
   \emph{Primes of the form $x^2 + ny^2$. 
  Fermat, class field theory and complex multiplication},
   a Wiley-Interscience Publication. John Wiley \& Sons, Inc., New York, 1989.
\bibitem{Daless}
  D'Alessandro, Flavio,
  \emph{Free groups of quaternions},  
  Internat. J. Algebra Comput. \textbf{14}(2004),  no.~1, 69--86.
\bibitem{Djokovic}
  Djokovi\'{c}, Dragomir \v{Z}.,
  \emph{A correction, a retraction, and addenda to my paper: 
  ``Another example of a finitely presented infinite simple group''}, 
  J. Algebra \textbf{82}(1983), no.~1, 285--293.
\bibitem{Dyubina}
  Dyubina, Anna,
  \emph{Instability of the virtual solvability 
  and the property of being virtually torsion-free for quasi-isometric groups},
  Internat. Math. Res. Notices \textbf{21}(2000), 1097--1101.
\bibitem{GAP}
  The GAP group, \textsf{GAP} --- Groups, Algorithms, and Programming, Version 4.4; 2004,
  \texttt{http://www.gap-system.org}
\bibitem{Johnson}
  Johnson, David L., 
  \emph{Presentations of groups}, 
  London Mathematical Society Student Texts, 15. 
  Cambridge University Press, Cambridge, 1990.
\bibitem{Mozes1}
   Mozes, Shahar,
   \emph{Actions of Cartan subgroups}, 
   Israel J. Math. \textbf{90}(1995), no.~1-3, 253--294. 
\bibitem{MR}
   Myasnikov, Alexei G.; Remeslennikov, Vladimir N., 
   \emph{Exponential groups. II. 
   Extensions of centralizers and tensor completion of CSA-groups},
   Internat. J. Algebra Comput. \textbf{6}(1996),  no.~6, 687--711. 
\bibitem{Rattaggi}
  Rattaggi, Diego,
  \emph{Computations in groups acting on a product of trees: 
  normal subgroup structures and quaternion lattices},
  Ph.D.~thesis, ETH Z\"urich, 2004.
\bibitem{Rattaggi2}
  Rattaggi, Diego,
  \emph{Anti-tori in square complex groups},
  Geom.\ Dedicata \textbf{114}(2005), 189--207.
\bibitem{Rattaggi3}
  Rattaggi, Diego,
  \emph{On direct product subgroups of $\mathrm{SO}_3(\mathbb{R})$},
  see arXiv:math.GR/0608292.
\bibitem{RaRo}
  Rattaggi, Diego; Robertson, Guyan,
  \emph{Abelian subgroup structure of square complex groups and arithmetic of quaternions},
  J. Algebra  \textbf{286}(2005),  no.~1, 57--68.
\bibitem{Stallings}
  Stallings, John R.,
  \emph{On torsion-free groups with infinitely many ends},
  Ann. of Math. (2) \textbf{88}(1968), 312--334.
\bibitem{Tits}
  Tits, Jacques,
  \emph{Free subgroups in linear groups},
  J. Algebra \textbf{20}(1972), 250--270.
\end{thebibliography}
\end{document}